\documentclass[a4paper]{article}
\usepackage{amsmath,amssymb}
\usepackage[ruled,vlined]{algorithm2e}
\usepackage{bm}
\usepackage{bbm} 
\usepackage[pdftex]{graphicx}
\usepackage{hyperref}
\usepackage{color}
\usepackage{cite} 

\usepackage{tikz}

\SetKw{KwBy}{by}

\definecolor{mywhite}{rgb}{1,1,1}

\newtheorem{theorem}{Theorem}

\newtheorem{lemma}[theorem]{Lemma}

\let\epsilon\varepsilon

\newcommand{\eps}{\epsilon}

\newcommand{\Luk}{\L{}ukasiewicz}

\newcommand{\cT}{\mathcal{T}}
\newcommand{\cY}{\mathcal{Y}}

\newcommand{\cO}{\mathcal{O}}

\newcommand{\cZ}{\mathcal{Z}}
\newcommand{\cX}{\mathcal{X}}

\newcommand{\cB}{\mathcal{B}}
\newcommand{\cM}{\mathcal{M}}

\newcommand{\bR}{\mathbb{R}}

\newcommand{\bI}{\mathbb{I}}

\newcommand{\bH}{\mathbb{H}}

\newcommand{\bZ}{\mathbb{Z}}

\newcommand{\bC}{\mathbb{C}}

\newcommand{\bN}{\mathbb{N}}

\newcommand{\bbu}{\mathbbm{1}}

\newcommand{\cyc}{\mathrm{Cyc}}

\newcommand{\bnu}{\bm{\nu}}

\newcommand{\bern}{\mathrm{Bern}}

\newcommand{\rint}{\mathrm{RndInt}}

\def\smfrac#1#2{{\textstyle\frac{#1}{#2}}}

\newcommand{\be}{\begin{equation}}
\newcommand{\ee}{\end{equation}}
\newcommand{\ef}[1]{\, #1}

\newcommand{\lam}{\lambda}

\newcommand{\dx}[1]{\!\mathrm{d}#1\,}
\newcommand{\dxns}[1]{\mathrm{d}#1}

\newcommand{\deenne}[2]{\frac{\partial^#2}{\partial #1 ^#2}}
\newcommand{\dienne}[2]{\frac{\mathrm{d}^#2}{\mathrm{d} #1 ^#2}}
\newcommand{\deennearg}[3]{\frac{\partial^#2 #3}{\partial #1 ^#2}}

\newcommand{\bt}{{\bm t}}

\newcommand{\eee}{\mathbb{E}}
\newcommand{\ppp}{\mathbb{P}}

\begin{document}





\begin{center}
{\Large \bf Boltzmann sampling of irreducible\\
\rule{0pt}{15pt}%
 context-free structures
in linear time}
\\
\rule{0pt}{20pt}%
{\large Andrea Sportiello}
\\
\rule{0pt}{14pt}%
LIPN, Universit\'e Sorbonne Paris Nord,
and CNRS, F-93430 Villetaneuse, France
\\
\rule{0pt}{13pt}%
{\tt
  andrea.sport\makebox[0pt][l]{i}iello@l\makebox[0pt][l]{i}ipn.un\makebox[0pt][l]{i}iv-paris13.fr}
\\
\rule{0pt}{14pt}%
\today
\end{center}

\vspace{5mm}

\noindent 
\rule{.075\textwidth}{0pt}%
\begin{minipage}{.85\textwidth}
{\small{\bf Abstract:} We continue our program of improving the
  complexity of so-called `Boltzmann sampling' algorithms, for the
  exact sampling of combinatorial structures, and reach average
  linear-time complexity, i.e.\ optimality up to a multiplicative
  constant.  Here we solve this problem for `irreducible context-free
  structures', a broad family of structures to which the celebrated
  Drmota--Lalley--Woods Theorem applies.
  Our algorithm is a `rejection algorithm'.  The main idea is to
  single out some degrees of freedom, i.e.\ write $p(x)=p_1(y)
  p_2(x|y)$, which allows to introduce a rejection factor at the level
  of the $y$ object, that is almost surely of order $1$.}
\end{minipage}
\vspace{4mm}



\section{The quest for linear-time Boltzmann sampling}
\label{sec.intro}

\noindent
We continue our study of the so-called ``Boltzmann''
Exact Sampling Algorithm \cite{flajBS}, a wide class of algorithms
which allow to sample random combinatorial structures of size $n$ from
a given measure, in an average time that scales algebraically with
$n$.  This algorithm has received many praises since its appearence in
2004, and it is fair to state that nowadays constitutes a small branch
by its own with the Theory of Algorithms. It exists in various
incarnations, such as labeled \cite{flajBS} and unlabeled
\cite{flajBS2} combinatorial structures, structures defined through
differential specifications \cite{bodi1ods} and multi-parametric
extensions \cite{bodiMulti}. In this paper we will mostly concentrate
on the `original' case, discussed in \cite{flajBS}, of labelled
structures.

The structures to which this family of algorithms mainly apply are
described in detail in the Flajolet and Sedgewick monograph on
Analytic Combinatorics \cite{flaj}, where a large stress is given to
generating-function techniques and saddle-point methods for the
asymptotic enumeration of the configurations.  In short, the Boltzmann
Algorithm explores the possibility of translating the informations
implicit in this analysis into an efficient exact-sampling
algorithm. Say that we have combinatorial objects
$x \in \cX = \bigcup_n \cX_n$, and objects of size $n$ have measure
$\mu_n(x)$, with support $\cX_n$. In most cases, a one-parameter
family of measures exist, $\mu(x;\beta)$, which takes the form
$\mu(x;\beta)=\sum_n p_{\beta}(n) \mu_n(x)$, and is `natural' in the
sense that it is the one implicit in the construction of the
generating function $X(z)$ for the objects (where $z$ and $\beta$ are
easily related). In other words, the measures $\mu_n(x)$ and
$\mu(x;\beta)$ are the canonical and grand-canonical Boltzmann--Gibbs
measures, at energy $n$ and inverse-temperature $\beta$, for the
statistical ensemble consisting of the combinatorial objects, and the
name `Boltzmann Method' is a tribute to the role of Ludwig Boltzmann
in the foundations of Statistical Mechanics, whose ideas are of
inspiration for the method.  The recursive description of $X(n)$
implicit in the combinatorial specification translates into a
linear-time algorithm for sampling from $\mu(x;\beta)$, which thus
induces an algorithm for sampling from $\mu_n(x)$, with complexity
$\sim n/p_{\beta}(n)$.

Essentially in all the cases of interest for us, the \emph{Shannon
entropy} of the measure $\mu_n(x)$, defined as
$S[\mu_n]:=-\sum_{x \in \cX_n} \mu_n(x) \ln \mu_n(x)$, scales linearly
with $n$, and, as well known, this provides a lower bound to the
complexity for exact sampling, 
because any algorithm needs to sample on average at least $S[\mu_n]$
random bits, and the cost of a random bit is of order 1. (Note that,
if the measure is the uniform measure, as in most of the concrete
applications, then $S[\mu_n]=\ln A_n$, and if $A_n \sim \rho^{-n}
n^{\gamma}$, then $S[\mu_n]=-n \ln \rho + o(n)$.)  We say that an
algorithm is
\emph{optimal up to a multiplicative constant}\footnote{This is a
  different, weaker notion w.r.t.\ optimality \emph{tout court},
  which, in the framework of exact-sampling algorithms, corresponds to
  optimality in the average number of random bits, up to subleading
  corrections, and a time complexity for other operations which is not
  larger than the complexity for sampling the required random
  bits. This version of optimality is the `ultimate goal' of exact
  sampling, but it is rather exceptional, and, in our opinion,
  presently beyond reach at the level of generality of the present
  paper. Among the few successful examples obtained so far we mention
  \cite{axel3, axel4, axel5, axel6} for various families of trees and
  walks, \cite{axel1, axel2, duchtree} for random permutations,
  \cite{bodi1} for ``linear extensions in series-parallel posets''
  among others.}
if it has average complexity bounded from
above by $C \cdot S[\mu_n]$, for some constant $C\geq
1$. Unfortunately, the Boltzmann algorithm \emph{as is} is not
optimal, as normally $\max_{\beta}p_\beta(n)$ scales as $1/n$, or, at
best, as $1/\sqrt{n}$. Thus, the average complexity scales as $\sim
n/p_{\beta}(n) = n^{1+\gamma}$ for some $\gamma>0$.  For this reason
we have started to explore under which circumstances we can improve
the ideas of Boltzmann sampling, up to (possibly) reach optimality.

In a first paper \cite{gascom16}, we devise a general method that applies
to all the cases in which an object $x$ of size $n$ can be decomposed
canonically as the union of two parts, $x=x_1 \cup x_2$, both of size
$\sim n$. When this is the case, an improvement of the na\"ive
Boltzmann method for sampling from the Hadamard product of two
distributions allows to decrease the complexity from $n^{1+\gamma}$ to
$n^{1+\frac{\gamma}{2}}$. In our opinion, this is already an
interesting result, as it is rather general, however it does not cover
all the applications of Boltzmann sampling (as the forementioned
canonical decomposition does not always exist), and, most importantly,
still does not reach optimality.

In a second set of works \cite{analco15,gascom18}, we describe how to exactly sample, in
linear time, combinatorial structures which (up to bijections) are
described by random bridges (i.e., random lattice walks from $(0,0)$
to $(n,0)$, with steps $(+1,h_j)$), even when the distributions
$p_{n,j}(h_j)$ of the steps are not all equal, in particular (as
described in
\cite{gascom18})
provided that they
satisfy a property that we call \emph{positive decomposition}. We show
that this technical constraint is not very restrictive, as the realm
of applications includes, among other things, important combinatorial
classes such as random set-partitions, that is partitions of a set of
$n$ elements into $k$ non-empty subsets, counted by Stirling numbers
of the second kind and related to the exact sampling of minimal
automata
\cite{BasDavSpoMA}, and its `dual' problem, of random permutations of
size $n$
with $k$ cycles, counted by Stirling numbers of the first kind, for
which, to our knowledge, no linear-time sampling algorithms were
previously available. However, admittedly, the problems solved by this
algorithm are just a relatively small subclass of all those that are
addressed by the Boltzmann method, so a quest for a new idea that
could tackle very large classes of combinatorial objects still stands.

Let us put the Boltzmann Method and its variants more in context. Of
course, this method
is not the only available exact-sampling algorithm for combinatorial
structures. One relatively ancient algorithm is due to Nijenhuis and
Wilf \cite{wilf} (with later modifications by Flajolet, Zimmermann and
Van Cutsem \cite{FlajZvC}, see also \cite{DenZimm}). If the exact
enumeration of certain combinatorial objects is described by a system
of equations, of the form
\[
A^{(i)}_n = \cdots + 
\!\!\!\!
\sum_{\substack{ 
n_1,n_2,\ldots,n_{\ell} \geq 0 \\ n_1 + n_2 + \cdots + n_{\ell} \leq
n}} 
\!\!\!\!
c(n_1,n_2,\ldots,n_{\ell},n) 
\,
A^{(j_1)}_{n_1} A^{(j_2)}_{n_2}
\cdots A^{(j_{\ell})}_{n_{\ell}} + \cdots
\]
where the coefficients $c$ are 
positive-valued functions 
(for their range of arguments), 
that can be evaluated as floats with
$\Theta(\ln n)$ digits in time $\Theta(\ln n)$,
then, for $(i_1,\ldots,i_{\ell})$ 
being sets of indices appearing as left factors of the monomials above,
the quantities
$A^{(i_1,\ldots,i_{\ell})}_n := \sum_{n_1,\ldots,n_{\ell}\,:\,\sum_j
  n_j=n} A^{(i_1)}_{n_1} \cdots A^{(i_{\ell})}_{n_{\ell}}$,
can be evaluated recursively up to the size of interest, once and for
all in a `preprocessing phase', in a time of order $n^2 \ln n$, and
the final data structure occupies a space of order $n \ln n$. At this
point, a simple divide-and-conquer algorithm (with complexity improved
by the so-called ``boustrophedon method'') allows to perform exact
sampling in average time of order
$n \ln n$. So, if one is interested in exact sampling with the aim of
performing a statistical average, when the number of samples is much
larger than the size of the objects (as is often the case), and allows
for extra logarithmic factors, even this ``brute-force counting''
strategy, once improved by the expedients outlined above, may be
quasi-optimal.

In fact, when the system of equations is linear (as in the broad
family of \emph{regular languages}), the setting of this algorithm
simplifies drastically, and a careful implementation of these ideas
allows to achieve linear complexity even for the single-run sampling,
that is, with an algorithm in which also the preprocessing is
quasi-linear \cite{Denise, BernGim}.\footnote{Note that in this case
  optimality is reached by carefully dealing with floating-point
  arithmetics, something which is doable and nowadays well-understood,
  however painful for what concerns algorithmic validation of the
  results, and, for ``purists'' of Complexity Theory, less elegant
  than a purely combinatorial algorithm. More precisely, we say that
  an algorithm `uses floats' if it requires the high-precision
  evaluation of quantities that also depend on the size, like for
  example, for a given random number
  $y \in [0,1]$, the largest integer $m$ such that $\sum_{j=0}^m
  2^{-n} \binom{n}{m}<y$, while we consider `legitimate' the
  evaluation, once and for all in the preprocessing phase, of a finite
  number of parameters, depending only on the combinatorial structure
  at hand, and not on the goal size, or also the use of potentially
  complicated quantities, such as $\binom{2n}{n} 2^{-2n} \sqrt{n}$, as
  acceptance rates in an algorithm, as in this case on average we need
  to certify only $\cO(1)$ binary digits of the quantity.}  So, we
should focus on combinatorial structures whose generating functions
are determined by \emph{non-linear} systems, which, indeed, are in
general more complicated classes (just like, in Algebra, non-linear
systems are more complicated than linear systems).  Within these
classes, there is an important family, in which the coefficients
$c(n_1,n_2,\ldots,n_{\ell},n)$ above are in fact constants, depending
only on $n-(n_1+n_2+\cdots+n_{\ell})$, and only finitely many are
non-zero. In this case one says that the corresponding structures are
\emph{context-free}. As explained in \cite[Sec.\;VII.1--6]{flaj},
under suitable hypotheses, the crucial \emph{Drmota--Lalley--Woods
  Theorem} (DLW)
\cite{Drmota, Woods, Lalley} applies, and these structures fall under
the so-called \emph{smooth inverse-function schema} (SIFS), that is
the generating function has a peculiar square-root singularity at the
critical point, and the enumeration has asymptotics of the form
$A_n \sim \rho^{-n} n^{-\frac{3}{2}} \exp(\cO(1))$, where the exponent
$-\frac{3}{2}$ is `universal', that is, it is shared by all classes
within the SIFS.

This family contains examples ranging from simple cases, such as
e.g.\ binary and unary-binary rooted planar trees, to rather
complicated ones, such as two-terminal planar graphs with no $W_5$
minor (a class of graphs included in the set of planar graphs, and
including series-parallel graphs).
Indeed, quoting \cite[pg.\;443]{flaj}:
\\
\rule{.02\textwidth}{0pt}%
\rule{0pt}{30pt}%
\begin{minipage}{.96\textwidth}
{\it there is a progression in the complexity of the schemas leading to
square-root singularity. From the analytic standpoint, this can be roughly rendered
by a chain}
[~inverse functions $\longrightarrow$ implicit functions
$\longrightarrow$ systems~].
{\it It is, however, often meaningful to treat each combinatorial problem at its minimal
level of generality}
\end{minipage}
\\ 
\rule{0pt}{14pt}%
This hierarchy of difficulty seems to stand also at the level of the
Boltzmann Method. In fact, the simplest case of the SIFS corresponds
to \emph{simple varieties of trees and inverse functions} in
\cite{flaj}, and is analysed in Section VII.3 there.  In this
situation, there exists an algorithm, due to Devroye \cite{devr}, that
achieves (quasi-)optimality,\footnote{I.e., it may have extra
  logarithmic factors, when bit complexity is taken into account.}  by
considering a classical bijection with \Luk\ 
paths, using the
cyclic lemma, and exploiting the exchangeability of the steps for the
corresponding bridges. This algorithm has the small flaw of involving
floating-point arithmetics, but it has the important merit of showing,
for the first time, that linear-time average complexity can be
achieved for `hard' (i.e., non-linear) combinatorial structures. Also,
when the step weights allow for a positive decomposition, our paper
\cite{gascom18} provides a different algorithm, that avoids floating-point
arithmetics (and the complicancies of dealing rigorously with it at
the level of programming).

However, there are two further steps in the complexity scale of the
SIFS, namely what is called
\emph{tree-like structures and implicit functions} in \cite{flaj}, and
treated in Section VII.4 there, and what is called 
\emph{irreducible context-free structures} in \cite{flaj}, and treated
in Section VII.6 there, which is mostly devoted to a discussion of the
forementioned Drmota--Lalley--Woods Theorem.  In this last case,
planar trees are replaced by a coloured variant (with as many colours
as equations in the system, plus one colour for the leaves), and the
natural bijection with lattice walks (i.e., the one induced by the
depth-first search countour of the tree) gives now walks with coloured
steps, and complicated non-local correlations between the various
steps. As a result, exchangeability is broken at the level of the
single steps, and the whole Devroye strategy cannot be applied,
without being complemented by some new idea. This fact is also
evidentiated in the original Devroye paper \cite{devr}, which explains
clearly to which situations his agorithm applies (and, implicitly, to
which situations his agorithm does \emph{not} apply).


The goal of this paper is to provide an average linear-time algorithm,
variant of the Boltzmann sampling method, that works in the broader
setting of irreducible context-free structures, thus performing ``two
leaps forward in one stroke'', on the complexity scale of the SIFS, in
the program of making the Boltzmann Method linear. Note however that
more complicated, non--context-free non-linear classes, such as what
is called
\emph{ordinary differential equations and systems} in \cite{flaj}, and
treated in Section VII.9 there, are still not covered by the treatment
of this paper, despite the fact that, as mentioned above, some linear
non--context-free problems are solved by our \cite{gascom18} (we hope to
address this level of generality in future work).

Our main idea is to decompose the combinatorial object $x$, in order
to extract one family of degrees of freedom which are specially
simple, and postpone the sampling of these degrees of freedom
\emph{after} the evaluation of the acceptance rate. By some magics
that we try to elucidate in Section \ref{sec.simpleex} on a simple
example (and discuss in full length in Section \ref{sec.anAR}), this
allows to gain the factor $1/p_{\beta}(n) \simeq n^{\gamma}$, and
reach optimality.  Ideas of this sort are not completely new (for
example, a version of this appears in our \cite{analco15,gascom18},
and another version appears in \cite{desalvo}), but are used here in a
different twist, that requires a number of subtle tweaks that we try
to introduce here in a pedagogical way.

The paper is organised as follows: besides the simple motivational
example of Section \ref{sec.simpleex}, in Section \ref{sec.ICFS} we
discuss some (more or less) well-known facts in Analytic
Combinatorics, concerning irreducible context-free structures. This
section is complemented by Section \ref{sec.ICFSexs}, that is mainly
devoted to examples.  Sections \ref{sec.cyclem} and \ref{sec.anaprel}
discuss some preliminary aspects of our algorithm (more of
combinatorial flavour in the first section, and of analytical flavour
in the second one). Finally, Section \ref{sec.algo} describes the
structure of the algorithm, and the functional form of the involved
quantities, while Section \ref{sec.anAR} describes how to optimise the
parameters, in such a way to reach a certification of optimal
complexity.


Section \ref{sec.anaprel} is complemented by an appendix that
discusses some facts in Perron--Frobenius Theory (which, as well
known, has a crucia role in the DLW Theorem). A second appendix
provides a reminder of the {\sc BalancedShuffle} algorithm
of Bacher, Bodini, Hollender and Lumbroso \cite{axel2}, which is used
as a black box within our algorithm. A third appendix collects some
technicalities required for the certification of optimality discussed
in Section \ref{sec.anAR}.

\section{A simple example}
\label{sec.simpleex}

\noindent
As a warm-up before introducing our full-fledged algorithm, let us
consider the exact sampling of lattice walks from $(0,0)$ to $(n,0)$,
with steps $(+1,\nu)$, with $\nu \in \{-1,0,+1\}$. Steps with $\nu=0$
come with a factor 2, that is, calling $x$ a walk, and $n_+$, $n_0$
and $n_-$ the number of steps $\nu=+1, 0, -1$ in $x$, respectively, we
have
\be
\mu_n(x) \propto 2^{n_0(x)}
\ef.
\ee
The weight $2$ has been chosen in order to have a trivial
normalisation: walks of this sort are just `walks of length $2n$ in
disguise': if $x'$ is a walk from $(0,0)$ to $(2n,0)$ with steps 
$\pm 1$, we obtain a map from walks $x'$ to $x$ by just taking steps
in pairs, and the factor $2^{n_0}$ is nothing but the number of
preimages under this map. As a result, we have the more explicit
\be
\mu_n(x) = \binom{2n}{n}^{-1} 2^{n_0(x)}
\ef,
\ee
and more generally
\be
A_{n,h}
:=
\sum_{x: (0,0) \to (n,h)} 2^{n_0(x)}
=
\binom{2n}{n+h}
\ef.
\ee
As a result, a simple divide-and-conquer algorithm allows to sample
these walks in linear time. We grow the walk step by step. Say that
the concatenation of the first $j$ steps has reached the position
$(j,h)$. Then we must continue with a step $\nu=+1, 0$ or $-1$ with
probabilities $A_{n-j-1,h+\nu}/A_{n-j,h}$,
that is, given by the triple
of rational functions
\[
\left(
\frac{(n-j-h)(n-j-h+1)}{(2n-2j)(2n-2j+1)}\, , \,
\frac{2 (n-j-h+1)(n-j+h+1)}{(2n-2j)(2n-2j+1)}\, , \,
\frac{(n-j+h)(n-j+h+1)}{(2n-2j)(2n-2j+1)}
\right)
\ef.
\]
Keeping probabilities
$(\frac{1}{4}, \frac{1}{2}, \frac{1}{4})$ as first approximation, and
refining the evaluation only if the sampling procedure requires it,
allows to avoid spurious logarithmic factors in bit complexity. 

So, in this case we have no need of inventing a new
algorithm. However, it is instructive to see other algorithms at work
here, where all the calculations can be performed explicitly, and a
number of subtleties are not required, before trying to generalising
their ideas to more complicated situations.

Before doing this, we shall remind that, for every $k$ of order 1, and
every $(n_1,\ldots,n_k)$ with $n=n_1+n_2+\cdots+n_k$, it is possible to sample uniformly random
shuffles $\sigma \in S_{n_1,\ldots,n_k}$ of the string
\[
(
\underbrace{1,1,\ldots,1}_{n_1},
\underbrace{2,2,\ldots,2}_{n_2},
\cdots ,
\underbrace{k,k,\ldots,k}_{n_k}
)
\]
(that is, permutations in 
$\mathfrak{S}_n / (\mathfrak{S}_{n_1} \times \cdots \times
\mathfrak{S}_{n_k})$) in a time linear in $n$.  If we allow for an
extra $\ln n$ factor in the complexity, this can be done just by
sampling a random permutation, e.g.\ with the classical Fisher--Yates
algorithm \cite{FY1,FY2}. Random-bit optimality (that is, using
$-\sum_{j=1}^k n_j \ln (n_j/n)+o(n)$ random bits) is reached by the
{\sc BalancedShuffle} algorithm of Bacher, Bodini, Hollender and
Lumbroso \cite{axel2}. We will need uniform random shuffles repeatedly
in the following, where it will be understood that `BBHL shuffling'
refers to this algorithm, and $\textrm{BBHL}(n_1,\ldots,n_k)$ is the
corresponding sampler.

\medskip

\noindent
{\bf Boltzmann Method.}
Let us see how the ordinary Boltzmann sampling
would work in this case. Call $\mu_{\infty}(\nu)$ the measure
$\mu_{\infty}(0)=\frac{1}{2}$, $\mu_{\infty}(\pm 1)=\frac{1}{4}$. We
just have:

\begin{algorithm}[H]
\Begin{
\Repeat{$\sum_j x_j=0$}
{
$x=(x_1,\ldots,x_n) \longleftarrow \mu_{\infty}^{\times n}$\;
}
\Return{$x$}
}
\caption{Boltzmann Method for $(\frac{1}{4}, \frac{1}{2}, \frac{1}{4})$-bridges.\label{algo.boltzex}}
\end{algorithm}

\noindent
Each run costs $\Theta(n)$ (and exactly $\frac{3}{2}n$ random bits on
average), and the probability that a run is accepted is
$4^{-n}\binom{2n}{n} \simeq 1/\sqrt{\pi n}$. So the overall complexity
is of order $n^{\frac{3}{2}}$, as anticipated.

\medskip

\noindent
{\bf Devroye Method.}  The method described in \cite{devr},
specialised to this case, works as follows. First, we sample
admissible triples $(n_+,n_0,n_-)$, with the appropriate probability
distribution. As we must have $\sum_j x_j=0$, admissible triples have
the form $(n_+,n_0,n_-)=(m,n-2m,m)$ for
$m\in \{0,1,\ldots,\lfloor n/2 \rfloor \}$, and the probability
distribution is
\be
p_n(m) = \binom{2n}{n}^{-1} \binom{n}{m\;m\;n-2m}\, 2^{n-2m}
\ee
Then, we perform a random shuffle of the string consisting of $m$
symbols $-1$, followed by $n-2m$ zeroes, and by $m$ $+1$'s:

\begin{algorithm}[H]
\Begin{
$m \longleftarrow p_n(m)$\;
$x=
(
\underbrace{-,\ldots,-}_{m},
\underbrace{0,\ldots,0}_{n-2m},
\underbrace{+,\ldots,+}_{m}
)
$\;
$\sigma \longleftarrow \textrm{BBHL}(m,n-2m,m)$\;
$x = \sigma \circ x$\;
\Return{$x$}
}
\caption{Devroye Method for $(\frac{1}{4}, \frac{1}{2}, \frac{1}{4})$-bridges.\label{algo.boltzex2}}
\end{algorithm}

\noindent
Sampling from $p_n(m)$ is complicated but feasible in sublinear time
(it is easily done in average time $\sim \sqrt{n}$, by sampling $y$
uniformly in $[0,1]$, one digit at the time as long as they are
needed, calculating once and for all $p_{n}(\lfloor n/4 \rfloor)$ as a
high-precision float, and summing up the $p_n(m)$'s, in order of
distance from $m=\lfloor n/4 \rfloor$, up to reach the threshold $y$,
which is done quickly as the ratio $p_n(m+1)/p_n(m)$ is a simple
rational function and w.h.p.\ we need $\cO(\sqrt{n})$ summands).  The
rest of the algorithm requires no other subtlety, and takes linear
time. So this algorithm is optimal.

\medskip

\noindent
{\bf Accelerated Boltzmann Method.}
Let us now see how our acceleration method improves the complexity of
Boltzmann sampling.  Call $\mu_{\infty}'(\nu)$ the measure
$\mu_{\infty}'(0)=\frac{2}{3}$, $\mu_{\infty}'(1)=\frac{1}{3}$. The
structure of our algorithm is as follows:

\begin{algorithm}[H]
\Begin{
\Repeat{$u=n$ and $\mathrm{rand}_{[0,1]}<r_n(n-t)$}
{
$t=1$\;
$u=0$\;
$x=(-1,-1,\ldots,-1)$\;
\While{$u<n$}
{$x_t \longleftarrow \mu_{\infty}'$\;
 $t \to t+1$\;
 $u \to u+1+x_t$\;
}
}
$\sigma \longleftarrow \textrm{BBHL}(t,n-t)$\;
$x = \sigma \circ x$\;
\Return{$x$}
}
\caption{Accelerated Boltzmann Method for $(\frac{1}{4}, \frac{1}{2}, \frac{1}{4})$-bridges.\label{algo.boltzex3}}
\end{algorithm}

\noindent
In words, we sample a random walk, with steps $0$ and $+1$, and
probabilities $2/3$ and $1/3$ (thus with average slope $1/3$), up to
reach or jump over the line passing through $(n,0)$ with slope
$-1$. The probability of jumping over is roughly $1/4$ (with
exponentially small corrections).\footnote{Because we jump over the
  $n$-th diagonal if we are in $n-1$, and we take a $+1$ step, so the
  probability $p$ of not occupying a diagonal satisfies the
  steady-state equation $p=(1-p)\frac{1}{3}$, that gives
  $p=\frac{1}{4}$.}  The most important fact is that, assuming that we
reached the line at position $(t,n-t)$, we keep what we have
constructed so far with a suitable acceptance rate $r_n(n-t)$, and
restart otherwise. Once this acceptance step has been verified, we
just sample a random shuffle, for shuffling the set of steps $+1$ and
$0$ altoghether with the set of steps $-1$, and we are done.

So, our algorithm is optimal up to a multiplicative constant if we can
determine a function $r_n(m)$ such that:
\begin{enumerate}
\item the algorithm samples according to the desired measure;
\item $r_n(m) \in [0,1]$ for all $n\geq 1$ and 
  $0 \leq m \leq \lfloor n/2 \rfloor$;
\item the function $r_n(m)$ can be calculated efficiently, i.e., in
  time at most $\cO(n)$;
\item $\eee(r_n(m))=\Theta(1)$.
\end{enumerate}
Let us first address the most obvious constraint: that we are sampling
the desired measure. We do this by analysing the probability of
getting any given output string $x$, with 
$(n_+,n_0,n_-)=(m,n-2m,m)$.
We have a factor $2^{n-2m}/3^{n-m}$ for sampling the $0/+1$ walk up to
the line, then a factor $\binom{n}{m}^{-1}$ for sampling the unique
shuffle that produces the string under investigation, and finally we
have the acceptance rate $r_n(m)$. The resulting
product must be proportional to $2^{n-2m}$. This gives the equation
\be
r_n(m) \frac{2^{n-2m}}{3^{n-m} \binom{n}{m}}
=
K_n \,2^{n-2m}
\ee
for some $K_n$. That is,
\be
r_n(m)
=
K_n \, 3^{n-m} \binom{n}{m}
\ef.
\ee
Then, we have to choose $K_n$ as large as possible (in order to have
good hopes on the fourth condition), while satisfying the second and
third condition, that is, we have to choose $K_n$ as large as
possible, while keeping it easy to evaluate, and certified to be
smaller than $1/\max_m (3^{n-m} \binom{n}{m})$. For all $n$, the
sequence $c_{n,m} = 3^{n-m} \binom{n}{m}$ is log-concave, so it has a
unique maximum, at the value $m^*$ where 
$c_{n,m^*}/c_{n,m^*-1} \geq 1$ and $c_{n,m^*+1}/c_{n,m^*} \leq 1$.
As we have
$c_{n,m+1}/c_{n,m} =\frac{n-m}{3(m+1)}$, this gives 
$m^*=\lfloor \frac{n-3}{4} \rfloor$. So we can choose
\be
r_n(m)
=
\frac{3^{m^*}\, (n-m^*)!\, m^*!}
{3^{m}\, (n-m)!\, m!}
\ef.
\ee
Typical values of $|m-m^*|$ are of order $\sqrt{n}$, thus calculating
the quantity above, as a $d$-digit float, by calculating the
corresponding Pochhammer functions,
takes on average $\Theta(\sqrt{n} (d+\ln n))$. Of
course, we can do much better.
Recalling the Robbins bound on factorials \cite{robbinsbound}
\be
\label{eq.robb}
e^{\frac{1}{12 n+1}}
\leq 
\frac{n!}{\sqrt{2 \pi} n^{n+\frac{1}{2}}e^{-n}}
\leq e^{\frac{1}{12 n}}
\ee
we can determine if a random uniform number in $[0,1]$ is larger or
smaller than the acceptance rate above, in time $\Theta(1)$, with
probability $1-\Theta(1/n)$, and then for the remaining probability we
can perform the product in the way outlined above, this giving overall
complexity
$\Theta(1) \times (1-\Theta(1/n)) + \Theta(\sqrt{n} \ln n) \times
\Theta(1/n) = \Theta(1)$.

So, we are ready to address the one and only subtle point, which shall
illustrate the reason of the acceleration. That is, we shall
understand why $\eee(r_n(m))=\Theta(1)$, while for ordinary Boltzmann
$\ppp(\sum_j x_j=0)=\Theta(n^{-\frac{1}{2}})$, just as a result of the
fact that we have sampled a 0/+1 random walk up to the line with
slope $-1$ passing through $(n,0)$, instead of sampling a $-$1/0/+1
random walk up to the vertical line passing through $(n,0)$.  

In fact, we would have had the very same complexity of ordinary
Boltzmann if we did perform the `wrong' na\"ive choice 
$r_n(m) = \frac{3^{n-m}}{4^n} \binom{n}{m}$ (which, by the way, is
manifestly bounded by 1 because $\sum_m r_{n}(m)=1$).  However, we
could push up the acceptance rate, by a factor that is the inverse of
the maximum over $m$ of the na\"ive function, and this quantity indeed
is of order $\sqrt{n}$.

For what concerns the evaluation of $\eee(r_n(m))$, we could just use
the CLT for showing that $(m-m^*)/\sqrt{n}$ is asymptotically normal
distributed, and then check that $r_n(m)$ is also asymptotically
Gaussian, scaled \emph{not} as to be normalised, but rather as to have
maximum value 1. That is, roughly,
\be
\label{eq.5634343987}
\eee(r_n(m)) \sim \int \frac{\dx{x}}{\sqrt{2 \pi n \sigma^2}} 
\exp(-x^2/(2 n \sigma^2))
\exp(-x^2/(2 n \tau^2))
=
\tau / \sqrt{\sigma^2+\tau^2}
\ef,
\ee
where $x=m-m^*$, the first Gaussian is the approximation of the
probability of reaching $m$, and the second Gaussian is the acceptance
rate. This CLT principle applies in general.  Furthermore, in our
simple example we can perform the calculations explicitly, as we have
\be
\eee(r_n(m))
=
\frac
{
\sum_m \binom{n-m}{m} \frac{2^{n-2m}}{3^{n-m}}
\frac{\binom{n}{m}}{\binom{n}{m^*}} 3^{m^*-m}
}
{
\sum_m \binom{n-m}{m} \frac{2^{n-2m}}{3^{n-m}}
}
\ee
The denominator is nothing but $(3+(-3)^{-n})/4 \simeq 3/4$, which,
incidentally, is also the probability of reaching the $n$-th diagonal. The
numerator is the slightly more complicated expression
$3^{m^*-n} (2n)! m^*! (n-m^*)!/n!^3$. Note that there are three
factorials in the numerator, and three in the denominator, so, as $m^*/n=\Theta(1)$, there are no
$\sqrt{2 \pi n}$ factors coming from Stirling approximation, and
indeed, for large $n$ the numerator converges to $\sqrt{3/8}$, with
corrections of order $1/n$. That is, combining numerator and
denominator,
$\eee(r_n(m)) \simeq \sqrt{2/3}$,
%
which in turns, multiplying by the probability $3/4$ of reaching the
line instead of jumping over (or, equivalently, omitting the
denominator), gives on average $\sqrt{8/3}$ tries for accepting a run
of the algorithm. This completes our proof, and (given the optimality
of BBHL shuffling, and of sampling i.i.d.\ random values from $\mu'_{\infty}$)
tells us that our
algorithm has average bit complexity $\sim \sqrt{8/3}\, S[\mu_n]$.

\section{Irreducible context-free structures and coloured trees}
\label{sec.ICFS}

\noindent
As explained in detail in the Analytic Combinatorics book \cite{flaj},
several interesting combinatorial structures admit a recursive
\emph{context-free} definition, that is, for every structure $y$,
weighted with the measure of choice $p(y)$, one can choose canonically
two integers $h$ and $k$, such that $y$ is decomposed into $h$
``atoms'' and $k$ sub-structures $y_1$,~\ldots,~$y_k$, with $p(y) =
\phi_{h,k} \prod_j p(y_j)$ and $|y|=h+\sum_j |y_j|$. In the symbolic
framework described in \cite{flaj} (generalised in the natural way for
dealing with weighted objects instead of just counting, see
e.g.\ \cite{duchon2011}), this reads
\be
\cY = \sum_{h,k} \phi_{h,k}\, \cZ^{\otimes h} \times \cY^{\otimes k}
\ef,
\ee
for $\phi_{h,k}$ real non-negative,
and, at the level of generating functions,
\be
Y(z) = \sum_{h,k} \phi_{h,k} z^h Y(z)^k
\ef.
\ee
The $\phi_{h,k}$'s must satisfy certain technical conditions, mostly of
summability (no fat tails), and non-triviality.  Calling
$\Phi(z,y)=\sum_{h,k} \phi_{h,k} z^h y^k$, the equation above reads
\be
\label{eq.7648763}
Y(z) = \Phi(z,Y(z))
\ef.
\ee
This is the situation called 
\emph{tree-like structures and implicit functions} in
\cite{flaj}.
Many concrete examples are of the form $\Phi(z,y)=z\, \phi(y)$, that
is, every decomposition involves a single extra elementary node. This
is the case of simply-generated (rooted planar) trees, when each node
counts as an unit, and of \Luk\ excursions \cite[p.\;74]{flaj}, that
is lattice paths in the upper-half plane with steps of the form
$(+1,h)$ for $h \geq -1$ (these two families are in simple bijection,
where the path describes the depth-first search countour of the tree).
This is the special case called 
\emph{simple varieties of trees and inverse functions} in \cite{flaj},
and, for what concerns exact sampling, if $\phi(y)$ is a polynomial,
is covered by Devroye algorithm, while more general paths or trees
(for example paths in the upper-half plane with steps in $(1,\pm 1)$
and $(3,0)$) are in the more general framework (in the example, with
$\Phi(z,y)=z y^2+z^3y+z$).

An even more general framework is one in which we have more than one
(but finitely many, say $m+1$) types of combinatorial structures
$\cY^{(\alpha)}$, with $\alpha=0,1,\ldots,m$. Again, for every structure
$y$ of type $\alpha$, weighted with the measure of choice, one can
choose canonically $m+2$ integers, $h$ and $k_0$, $k_1$, \ldots, $k_m$, such
that $y \in \cY^{(\alpha)}$ is decomposed into $h$ ``atoms'' and
$k_{\beta}$ structures $y^{(\beta)}_1$, \ldots,
$y^{(\beta)}_{k_\beta} \in \cY^{(\beta)}$, for all $\beta=0,1,\ldots,m$, with
$p_{\alpha}(y) = \phi^{(\alpha)}_{h;k_0,k_1,\ldots,k_m} \prod_{\beta}
\prod_{j=1}^{k_{\beta}} p_\beta(y^{(\beta)}_j) $
and
$|y|=h+\sum_\beta \sum_j | y^{(\beta)}_j |$. 
In the symbolic framework described in
\cite{flaj}, this reads
\be
\cY^{(\alpha)} = \sum_{h,k_0,k_1,\ldots,k_m} 
\phi^{(\alpha)}_{h;k_0,k_1,\ldots,k_m}
\, \cZ^{\otimes h} 
\times (\cY^{(0)})^{\otimes k_0}
\times (\cY^{(1)})^{\otimes k_1}
\times
\cdots
\times
(\cY^{(m)})^{\otimes k_m}
\ef,
\ee
and, at the level of generating functions,
\be
Y^{(\alpha)}(z) = \sum_{h,k_0,k_1,\ldots,k_m} 
\phi^{(\alpha)}_{h;k_0,k_1,\ldots,k_m}
z^h \prod_{\beta} 
Y^{(\beta)}(z)^{k_\beta}
\ef.
\ee
Calling $\vec{Y}=(Y^{(0)},Y^{(1)},\ldots,Y^{(m)})$
and
$\vec{\Phi}=(\Phi^{(0)},\Phi^{(1)},\ldots,\Phi^{(m)})$,
and introducing the functions
$\Phi^{(\alpha)}(z,\vec{y})=\sum_{h,k_0,k_1,\ldots,k_m} 
\phi^{(\alpha)}_{h;k_0,k_1,\ldots,k_m}
z^h \prod_{\beta} 
y_{\beta}^{k_\beta}$, 
the equation above is just the natural multi-component version of
(\ref{eq.7648763}), that is
\be
\label{eq.2664576743}
\vec{Y}(z) = \vec{\Phi}(z,\vec{Y}(z))
\ef.
\ee
This is the situation called 
\emph{context-free structures and polynomial systems} in \cite{flaj},
and the main object of interest in this paper.  By convention (see
e.g.\ \cite[ex.\;I.54, pg.\;82]{flaj}), the combinatorial class to
which we are interested is the one represented by the first component
of our vector.

Just like equations of the form $Y=z \,\phi(Y)$ describe trees counted
with the number of nodes, and more generally equations of the form
$Y=\Phi(z,Y)$ can be related to trees where internal nodes and leaves
are distinguished, counted with the number of
leaves,\footnote{In order to have finitely many configurations for
  each given size, we require that no unary node can have an internal
  child, i.e., that $\Phi(z,Y)$ has no monomial $z^0 Y^1$.}
configurations associated to a context-free structure can be put in
bijection with certain `weighted coloured trees', in which the size is
the number of leaves, and the internal nodes can be `coloured' with
the indices from $0$ to $m$, and an internal node of colour $\alpha$,
with $k_{\beta}$ children of colour $\beta$ and $h$ children leaves
comes with a factor $\phi^{(\alpha)}_{h;k_0,k_1,\ldots,k_m}$ in the
weight.\footnote{Now, in order to have finitely many configurations
  for each given size, we require that the linear part of
  $\vec{\Phi}(0,\vec{Y})$ is a nilpotent matrix.}

\begin{figure}[tb]
\setlength{\unitlength}{85pt}
\begin{picture}(5,1.9)
\put(0,.9){\includegraphics[scale=.55]{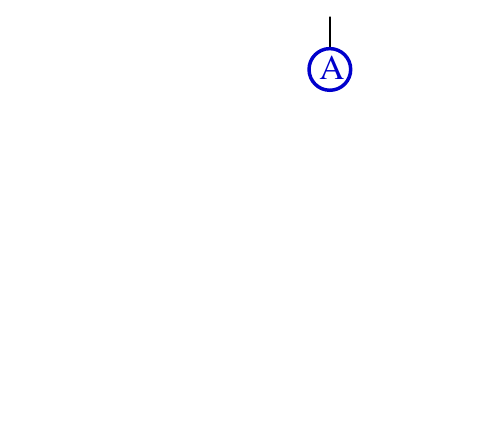}}
\put(1,.9){\includegraphics[scale=.55]{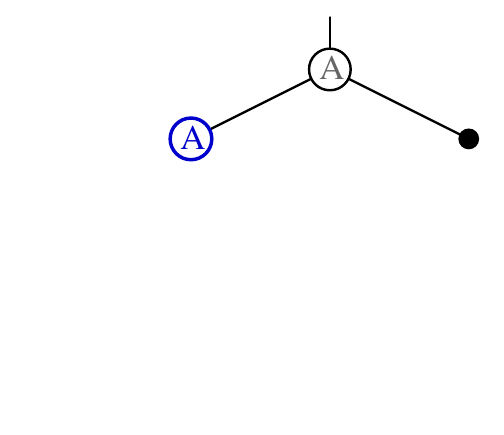}}
\put(2,.9){\includegraphics[scale=.55]{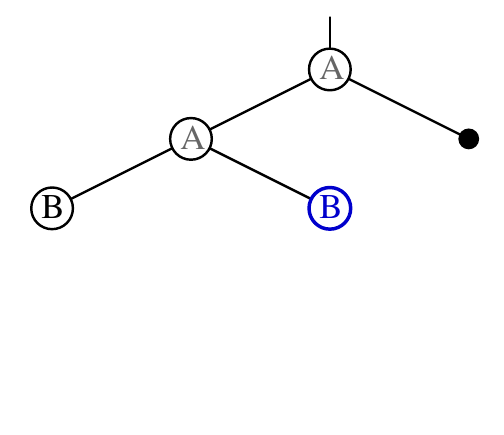}}
\put(3,.9){\includegraphics[scale=.55]{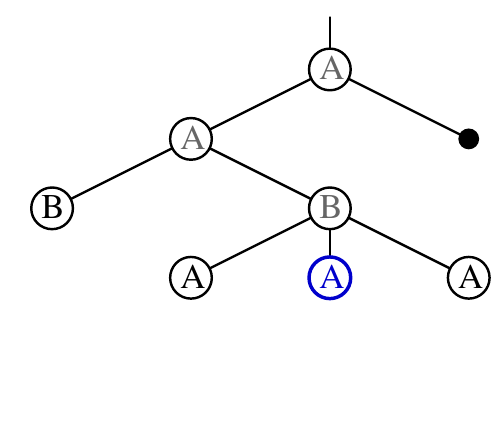}}
\put(4,.9){\includegraphics[scale=.55]{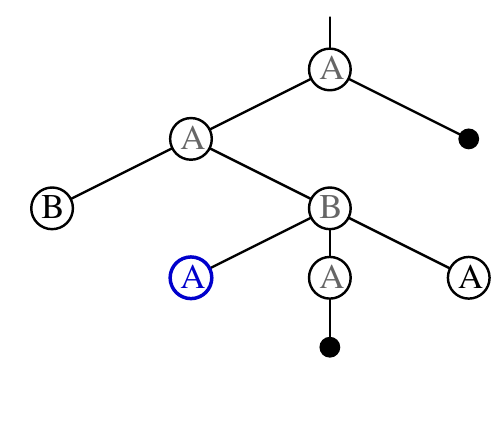}}
\put(0,0){\includegraphics[scale=.55]{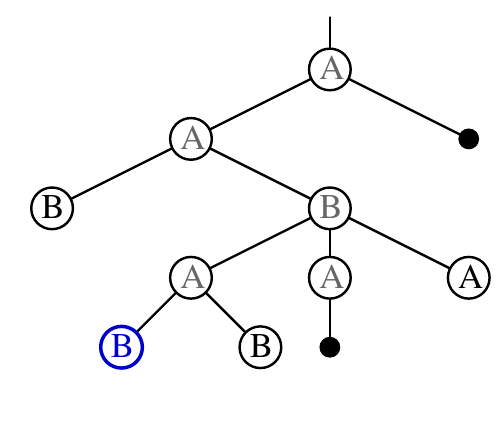}}
\put(1,0){\includegraphics[scale=.55]{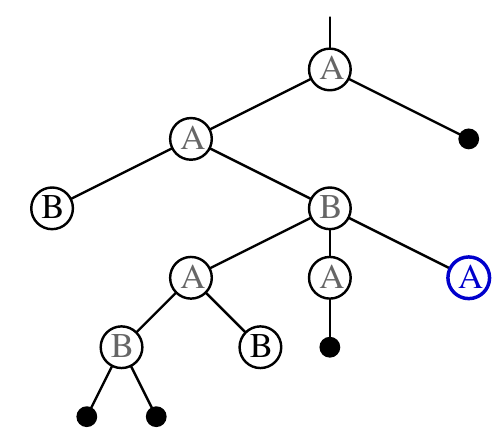}}
\put(2,0){\includegraphics[scale=.55]{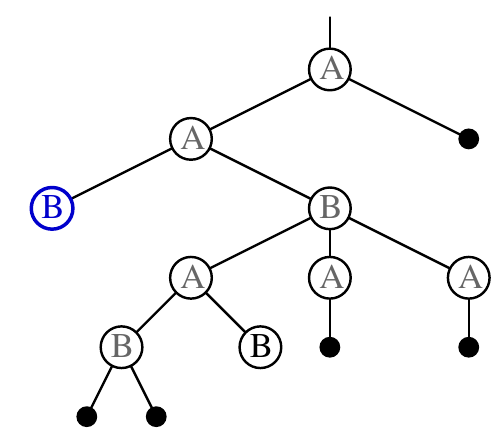}}
\put(3,0){\includegraphics[scale=.55]{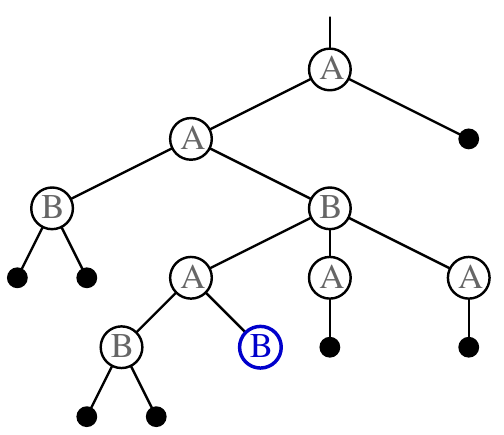}}
\put(4,0){\includegraphics[scale=.55]{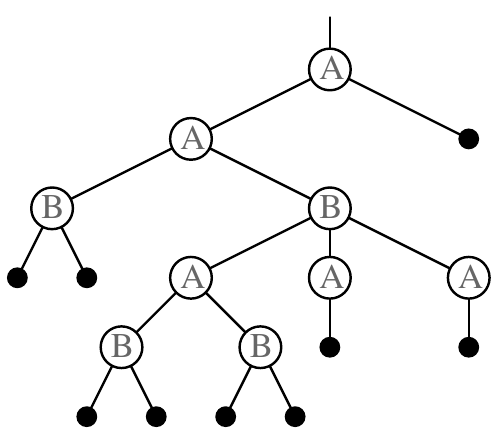}}
\end{picture}
\caption{\label{fig.rewrit}
Example of trajectory for the stochastic rewriting rules associated to
the specification
$A = A\,z+B^2+z\,;$ $B = A^3+z^2$.
In blue, the node in the stack that is being processed at the present
step of the process; in black, the other nodes in the stack; in gray,
nodes that have already been processed.}
\end{figure}

The specification given by the system can be interpreted as a
stochastic rewriting rule, of which every trajectory can be translated
into a coloured tree (see Figure \ref{fig.rewrit}). The process is
parametrised by a solution $(\rho,\vec{\tau})$ of the system
$\vec{\Phi}(\rho,\vec{\tau})=\vec{\tau}$.  It starts with a node of
label 0 at the root, which is the only node in a stack of `boundary
nodes'. Then, for a node in the stack with label $\alpha$, we choose
the composition of its offspring according to the probability
distribution
\be
p_{\alpha}(h;k_0,k_1,\ldots,k_m)
=
\tau_{\alpha}^{-1}
\phi^{(\alpha)}_{h;k_0,k_1,\ldots,k_m} \rho^h \prod_{\beta}
\tau_{\beta}^{k_{\beta}}
\ef.
\ee
(Note that, indeed, these probabilities are normalised.)  Then, the
node leaves the stack, and all of its $k_0+\cdots+k_m$ non-leaf
descendents are put in the stack. We continue the process, picking up
nodes from the stack, e.g.\ in random order or in a breath-first
search. The process stops when the stack is empty (or doesn't stop at
all).  The order of the operations does not affect the probability
distribution of the outcome, as long as it is guaranteed that, for
each finite height $h$, almost surely every node in the stack, at
height $h$ in the tree, is processed at some point, and in particular
this is the case when the process stops almost surely. The resulting
process, besides the minor complicancies coming from the colouring, is
essentially a Galton--Watson (GW) process. In particular,
it is well known that we have a critical Galton--Watson process
whenever these parameters correspond to the solution of the
\emph{characteristic system} \cite[pg.\;483]{flaj}, i.e., the solution
to the system of $m+2$ equations
\begin{align}
\label{eq.charsys}
\vec{\Phi}(\rho,\vec{\tau})&=\vec{\tau}
&
\det \left(
\bI - K
\right)
&=0
&
K_{\alpha \beta}
&:=
\deenne{\tau_{\beta}}{{}}
\Phi^{(\alpha)}(\rho,\vec{\tau})
\end{align}
in $(\bR^+)^{m+2}$ with smallest value of $\rho$,
while we have a subcritical GW process whenever
the spectrum of the matrix $K$ is strictly contained in the disk of
radius 1 (or, equivalently, the Frobenius eigenvalue of the matrix $K$
is strictly smaller than $1$). Indeed, the Frobenius eigenvalue of $K$
corresponds to the average number of children, in the GW process in
which the nodes of the stack are taken randomly.  Subcritical GW
processes lead to a probability distribution on the extensive
parameters of the tree (such as the number of nodes of a given type)
which has exponential tails, while, if the combinatorial system is
also ``irreducible'', critical GW processes lead to a probability
distribution with tail
$p(n) \sim C \cdot n^{-\frac{3}{2}}$, with $C$ determined by the DLW
Theorem (the precise notion of irreducibility is presented in the
context of this theorem, we refer to \cite{flaj} for the details).
Supercritical GW processes, associated to the case in which the
Frobenius eigenvalue is larger than 1, lead to trees which have a
finite probability of being of infinite size. In this case, the
measure described by the parameters $(\rho,\vec{\tau})$, and
conditioned to produce finite trees, is well-defined, and interesting
in several respects, however we do not need this notion in this paper,
so we do not discuss further the supercritical case.

\section{Examples of irreducible context-free structures}
\label{sec.ICFSexs}

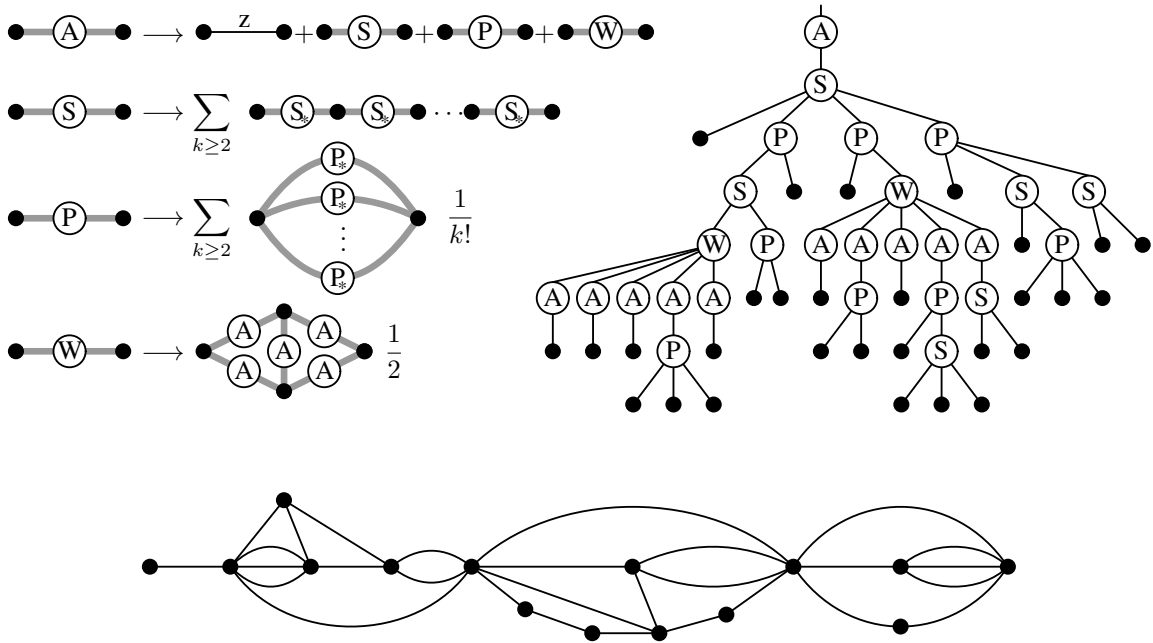
\begin{figure}[tb]
\setlength{\unitlength}{30pt}
\begin{picture}(10,6)(-.4,0)
\put(-.7,0){\includegraphics[scale=1]{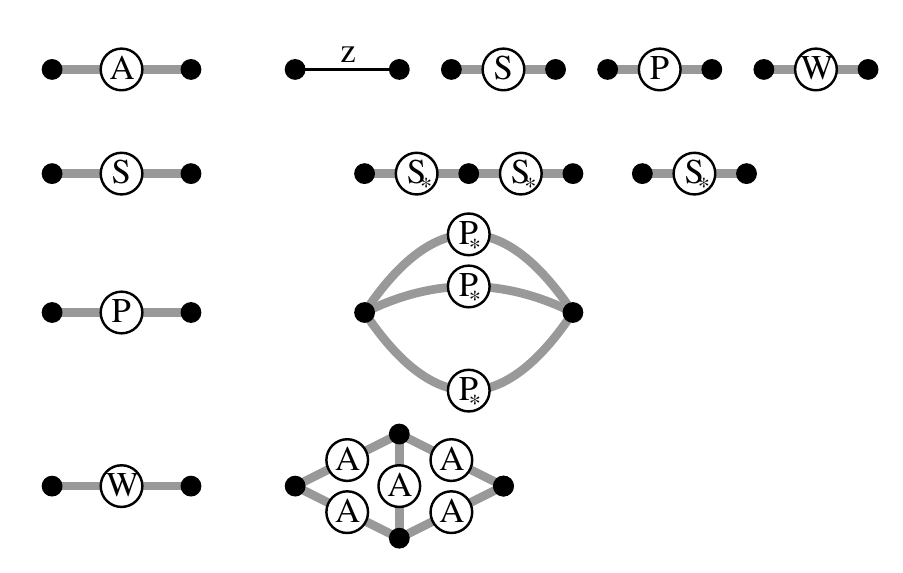}}
\put(6,0){\includegraphics[scale=1]{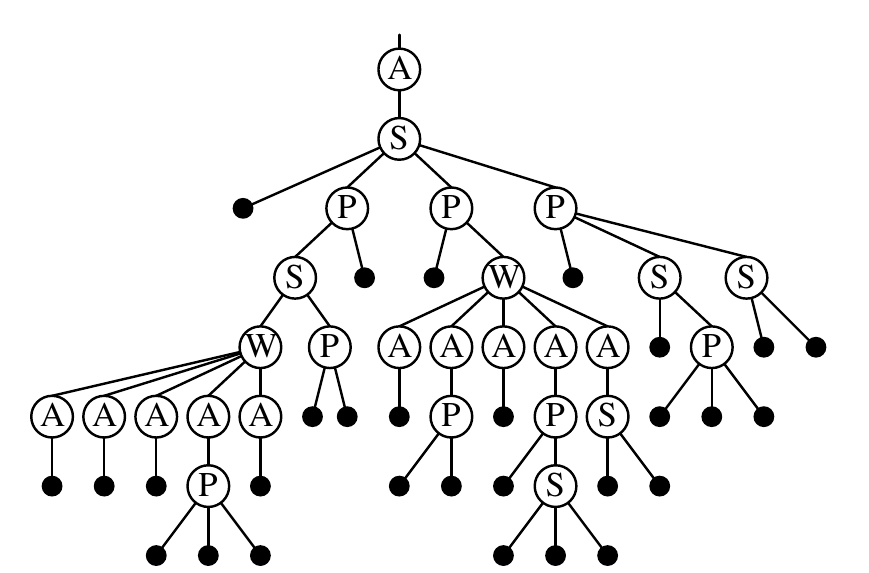}}
\put(1.37,.73){$\displaystyle{\longrightarrow 
\rule{72pt}{0pt}
\frac{1}{2}}$}
\put(1.37,2.43){$\displaystyle{\longrightarrow \sum_{k \geq 2}
\rule{39pt}{0pt} \raisebox{-10pt}{\vdots}
\rule{37pt}{0pt}
\frac{1}{k!}}$}
\put(1.37,3.75){$\displaystyle{\longrightarrow \sum_{k \geq 2}
\rule{74pt}{0pt} \raisebox{0pt}{$\cdots$}}$}
\put(1.37,4.73){$\longrightarrow$}
\put(3.27,4.74){$+$}
\put(4.77,4.74){$+$}
\put(6.27,4.74){$+$}
\end{picture}
\[
\includegraphics[scale=1]{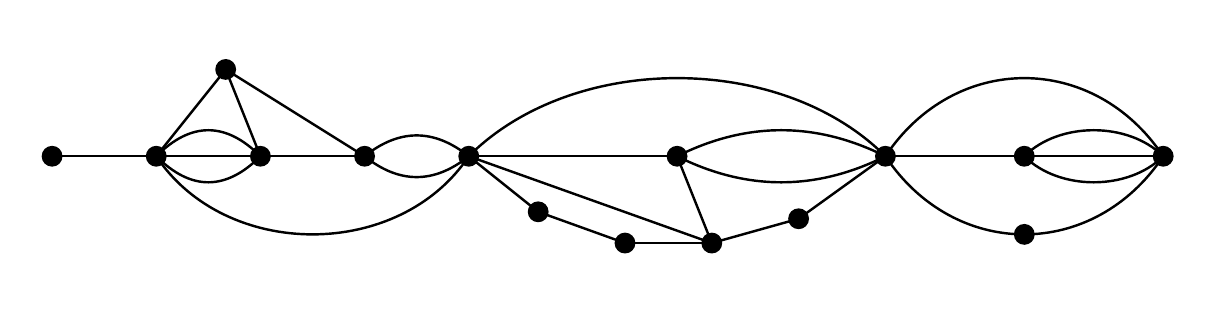}
\]
\caption{\label{fig.SPWex}Top left: the rewriting rules for 
two-terminal graphs with
no $W_5$ minor. As said in the text,
letters $A$, $P$, $S$ and $W$ 
denote
``all'', ``series'', ``parallel'' and ``Wheatstone bridge''
subclasses. Letters $S_{\ast}$ and $P_{\ast}$ denote 
$A \smallsetminus S=z+P+W$ and 
$A \smallsetminus P=z+S+W$, respectively.
Bottom: an example of configuration.
The two terminals are the left-most and right-most vertices.
Top right: the decomposition tree associated to the example.
}
\end{figure}

\noindent
A nice example of context-free structure is the class of
(two-terminal) series-parallel graphs, described for example in
\cite[p.\;72, ex.\;I.46]{flaj}, and more in detail in
\cite{PivSalSor2012}. Another, slightly more complex example (but also
more ``typical'', as, contrarily to series-parallel graphs, does not
have colourings alternating along the layers of the tree, and is
irreducible and aperiodic), is the class of (two-terminal) graphs with
no $W_5$ minor.\footnote{$W_5$ is the wheel graph with five vertices.}
To our knowledge, this class has been discussed so far only in a
seminar of ours, on the very same topic of this paper.\footnote{See
  {\tt
    https://library.cirm-math.fr/Record.htm?idlist=2\&record=19286312124910045949},
recorded during the meeting \emph{AofA: Probabilistic, Combinatorial
  and Asymptotic Methods for the Analysis of Algorithms}, on June
24th, 2019, at the \emph{Centre International de Rencontres
  Math\'ematiques} (CIRM), Marseille, France. Slides are available at
{\tt https://www.cirm-math.fr/RepOrga/1940/Slides/Sportiello.pdf}, and
the system (\ref{eq.3765746}) is on page 65. Note that in this
document there is a typo ($W_4$ in place of $W_5$) in the
excluded-minor description of the class of graphs.}
This class is described by the system of equations
\be
\label{eq.3765746}
\left\{
\begin{array}{l}
A = z + S + P + W \\
\rule{0pt}{10pt}%
S = \sum_{k \geq 2} (z + P + W)^k \\
\rule{0pt}{10pt}%
P = \sum_{k \geq 2} \frac{1}{k!} (z + S + W)^k \\
\rule{0pt}{10pt}%
W = \frac{1}{2} A^5
\end{array}
\right.
\ee
where letters $A$, $P$, $S$ and $W$ have been chosen to denote
``all'', ``series'', ``parallel'' and ``Wheatstone bridge'' subclasses
of these graphs.  See Figure \ref{fig.SPWex} for a description of the
specification as a rewriting system (left), a typical example of
configuration (bottom), and its representation as a coloured tree
(right).

Strictly speaking, these classes are not context-free structures,
because the functions $\Phi^{(\alpha)}$ are not polynomial.
Nonetheless, they can be treated on a similar ground, as in fact the
DLW Theorem has several extensions (see the discussion in
\cite[pg.\;493]{flaj}, in particular the Remark VII.29).

\begin{figure}[tb]
\[
\setlength{\unitlength}{50pt}
\begin{picture}(5,3)
\put(0,0){\includegraphics[scale=1]{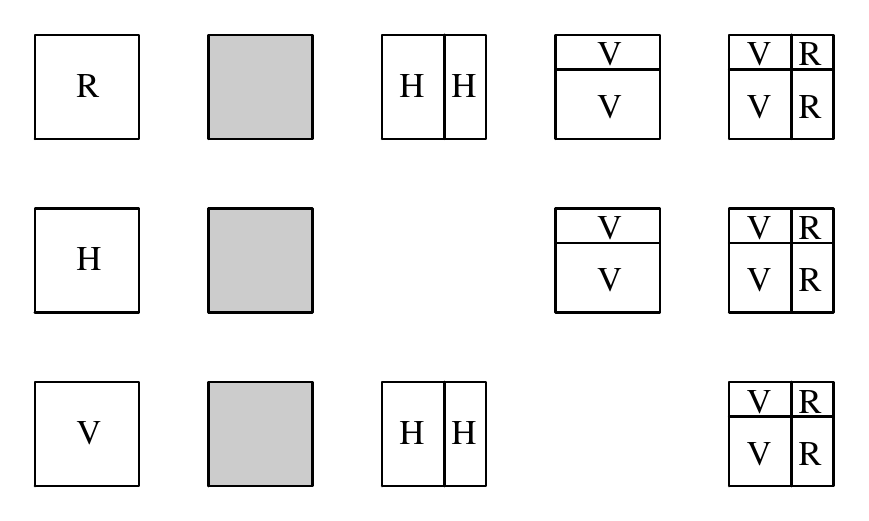}}
\put(0.92,2.45){$=$}
\put(0.92,1.45){$=$}
\put(0.92,0.45){$=$}
\put(1.94,2.45){$+$}
\put(1.94,0.45){$+$}
\put(2.94,2.45){$+$}
\put(2.94,1.45){$+$}
\put(3.94,2.45){$+$}
\put(3.94,1.45){$+$}
\put(3.94,0.45){$+$}
\end{picture}
\]
\caption{\label{fig.RHVrules}Graphical illustration of the
  combinatorial specification of `recursive topological
  subdivisions', as encoded by the system (\ref{eq.RHVsys}).}
\end{figure}
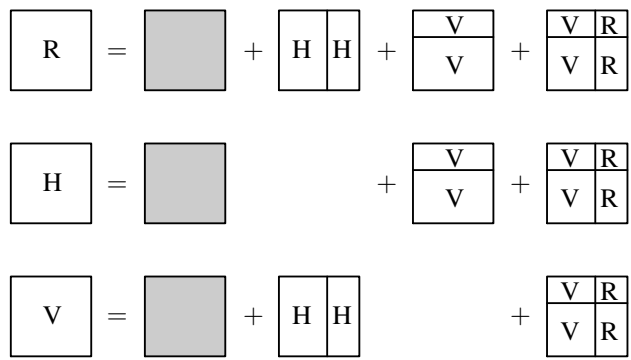

\begin{figure}[tb!]
\[
\includegraphics[width=\textwidth]{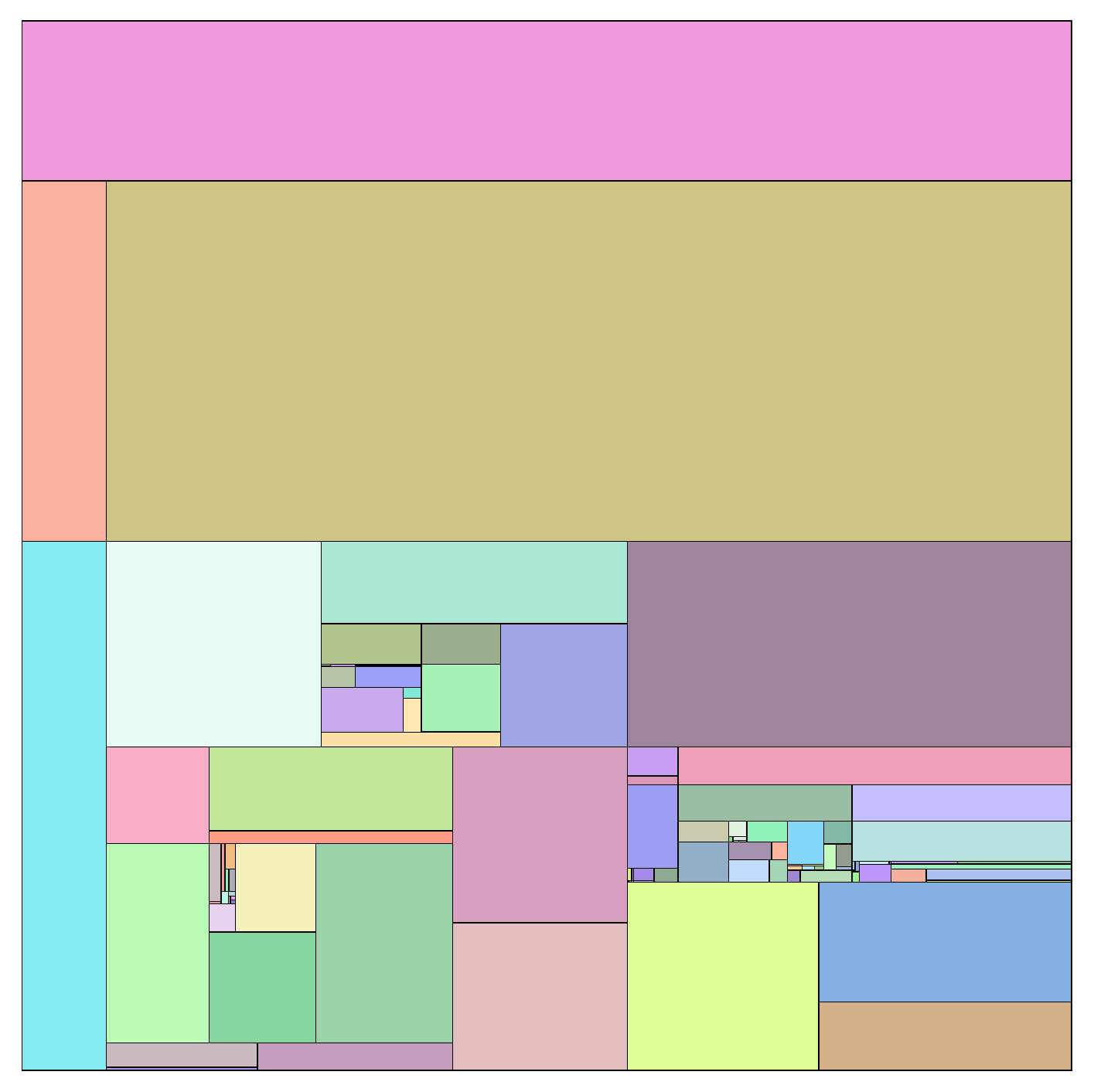}
\]
\caption{\label{fig.RHVex}Typical example of `recursive topological
  subdivision' of the unit square, of size 100.}
\end{figure}

Another example is constituted by recursive topological subdivisions
of a rectangle by straight lines. This is a new model, somewhat
pictorially similar to `quadtrees' 
(cf.\ \cite[Ex.\;VII.23, pg.\;523]{flaj}), and, if we trade rectangles
for triangles, to `stack-triangulations' \cite{AlbeMarck}.  In this
problem, every rectangle can be subdivided either horizontally, or
vertically, or in both directions. However we impose a further
condition, that in a sense avoids multiple countings of the same
configuration: if we have divided a rectangle horizontally, we cannot
divide horizontally any of the two resulting rectangles, and similarly
for vertical subdivisions.  See a typical example in Figure
\ref{fig.RHVex}. The size of a configuration is the number of
rectangles in the subdivision.  If we call $R(z)$ the generating
function for these configurations, and $H(z)$ (resp.\ $V(z)$) the ones
for rectangles that come from a horizontal (resp.\ vertical)
subdivision, these functions satisfy the system of equations
\be
\label{eq.RHVsys}
\left\{
\begin{array}{l}
R = z + H^2 + V^2 + R^4 \\
\rule{0pt}{10pt}%
H = z + V^2 + R^4 \\
\rule{0pt}{10pt}%
V = z + H^2 + R^4 
\end{array}
\right.
\ee
which is illustrated by Figure \ref{fig.RHVrules}.
If we use
the symmetry induced by (say) 90-degree rotations, we can identify the
generating functions $H$ and $V$, and reduce the system to two
equations
\be
\left\{
\begin{array}{l}
R = z + 2 H^2 + R^4 \\
\rule{0pt}{10pt}%
H = z + H^2 + R^4
\end{array}
\right.
\ee
We can further eliminate $H$, to get the algebraic relation
\be
1+4R = (1+z+R+R^4)^2
\ef,
\ee
that is, the inverse of the generating function, $z(R)$, is given by
\be
z(R)=\sqrt{1+4R}-1-R-R^4
\ef.
\ee
The critical values $(z^*,R^*,H^*,V^*)$ are thus given as the roots of
integer-valued polynomials (e.g., $R$ is determined by the equation
$z'(R)=0$), identified by the constraint of being real-positive, and
$z$ being the smallest real-positive root:
\begin{align}
\begin{split}
z^*:&\  0.1868943725\ldots \qquad 
16777216\, z^7
+113311744\, z^6
+276168704\, z^5
+277020672\, z^4
\\
&\  \phantom{0.1868943725\ldots} \qquad \quad 
+166353408\, z^3
-146027008\, z^2
+69783040\, z
-9433303 = 0
\end{split}
\\
R^*:& \ 0.3945155166\ldots \qquad 
64 R^7+16 R^6+32 R^4+8 R^3+4 R-3 = 0
\\
H^*:&\ 0.3028172374\ldots \qquad 
8 H^7+28 H^6+36 H^5+20 H^4+4 H^3+2 H-1=0
\qquad (V^*=H^*)
\end{align}

\section{From trees to walks, and the cyclic lemma}
\label{sec.cyclem}

\noindent
From this point onward, we assume that our context-free structure is
irreducible and aperiodic in the sense of the DLW Theorem. We also use
$A(z)$ as a synonim of $Y^{(0)}(z)$, as this first component, besides
being associated to the combinatorial class that we want to sample, plays
a distinct role in the whole construction of the algorithm.

As we described above, the combinatorial structures we shall sample
are in bijection with rooted planar trees $T \in \cT_n$, whose root
index is $A$, and which have $n$ leaves. The number of internal
nodes is not fixed, however it is bounded by $Kn$ for some constant
$K$ (as we have required that the linear part of $\vec{\Phi}$ is a
nilpotent matrix).

Several algorithms for the exact sampling of trees, including the
Devroye Algorithm among various others, make use of a bijection with
suitable lattice paths, and the cyclic lemma. However, as we said
above, this strategy does not apply immediately to coloured trees.

A better idea is to decompose the trees into subtrees, by breaking it
at all the internal nodes with label $A$. Each tree is thus described
by a list of subtrees, $T \equiv \bt = (t_1,\ldots,t_k)$, where a
subtree $t_j$ has root index $A$, leaves of two types ($A$ and $z$),
and all other internal nodes with labels in $\{1,\ldots,m\}$. We call
$\cT_{v,\ell}$ the class of such subtrees having $v$ leaves with label
$z$, and $\ell$ leaves with label $A$, and 
$\cT^{\rm (sub)}= \bigcup_{v,\ell \geq 0} \cT_{v,\ell}$ the class off all
subtrees altogether. Let us adopt the shortcuts
$v_i=v(t_i)$ and $\ell_i=\ell(v_i)$, and $k=k(\bt)$. A set of
necessary and sufficient conditions for a list $\bt$ to determine a
tree is the following:
\begin{enumerate}
\item
$
\sum_{j=1}^k v_j = n$;
\item
$
\sum_{j=1}^k (\ell_j-1) = -1$;
\item
$\sum_{j=1}^{k'} (\ell_j-1) \geq 0$ for all $k'<k$.
\end{enumerate}
In other words, the concatenation of the vectors
$\vec{u}_j:=(v_j,\ell_j-1)$ must constitute a generalisation of the
notion of \Luk\ path, in which the steps can go forward by an amount
different from $+1$. In lack of a name that has already been fixed in
the literature, we will call them \emph{generalised \Luk\ excursions}.

This construction has traded the complicancy of coloured steps, and
non-local correlations, with the (comparatively minor) complicancy of
steps with variable horizontal length. In particular, if we remove
the constraint of remaining in the upper-half plane, the corresponding
generalised \Luk\ walks have exchangeable steps.

If we call $A_n$ the partition function at size $n$, i.e.\ $A_n=[z^n]A(z)$,
we have
\begin{align}
\mu_n(T) &= \frac{W(T)}{A_n}\, \bbu(T \in \cT_n) 
&
W(T) &=
\prod_{\substack{
v \in V(T) \\
\textrm{internal}}}
\phi^{(\alpha(v))}_{h(v);k_0(v),k_1(v),\ldots,k_m(v)}
\end{align}
and, rephrased in terms of lists of subtrees
\be
\label{eq.37645734}
\mu_n(\bt) = \frac{\prod_j W(t_j)}{A_n}\, \bbu(
\textrm{$\bt$ satisfies the conditions (1), (2) and (3) above})
\ef.
\ee
At this point, we recall the Cyclic Lemma
\begin{lemma}[Cyclic Lemma]
Let $U$ be a set of steps in $\bZ^2 \smallsetminus (0,0)$, strictly contained in a
half-plane 
$\bH_{a/b}=\{(x,y)\,|\,ya+xb>0\}$.
Let $w_0=(\vec{u}_1,\ldots,\vec{u}_k) \in U^k$ be a walk from 
the origin to the point
$(n_x,n_y) \in \bH_{a/b}$.
%
For $0 \leq j < k$, call
$w_j=(\vec{u}_{j+1},\ldots,\vec{u}_k,\vec{u}_1,\ldots,\vec{u}_j)$ the
cyclic shifts of $w_0$.  If $n_x$ and $n_y$ are coprime, there exists
exactly one index $0 \leq j < k$ such that the walk $w_j$ is contained
in the half-plane $\bH_{-n_y/n_x}$.
\end{lemma}
We shall call $\cyc$ the operator such that, applied to a walk $w_0$,
gives the only cyclic shift $w_j$ of $w_0$ with the property above. An
immediate and crucial consequence of this lemma (also based on the
fact that our generalised \Luk\ excursions reach the point $(n,-1)$ and all
integers are coprime with $-1$), is that we can rewrite equation
(\ref{eq.37645734}) getting rid of the complicated condition (3):
\be
\label{eq.37645734b}
\mu_n(\bt) = 
\sum_{\bt' \,:\, \cyc \,\circ\, \bt'=\bt}
\frac{1}{A_n}\, 
\frac{\prod_j W(t'_j)}{k(\bt')}\,
\bbu(
\textrm{$\bt'$ satisfies the conditions (1) and (2) above})
\ef.
\ee
We will call \emph{generalised \Luk\ bridges} the corresponding
family of directed lattice walks, that is, the same walks as the
excursions, without the constraint of remaining in the upper half-plane.
As a consequence, if we are able to exactly sample lists $\bt'$ with
the measure\footnote{Here $\delta(n)$ denotes the Kronecker delta
  $\delta_{n 0}$, and \emph{not} the Dirac delta. This choice is based
  on the fact that here $n$ is replaced by large expressions, which would be
  poorly rendered as indices, and also there is no risk of confusions,
  as there are no Dirac delta's in this paper.}
\be
\label{eq.37645734c}
\begin{split}
\mu_n^{\rm (cyc)}(\bt') &= 
\frac{1}{A_n}\, 
\frac{\prod_j W(t'_j)}{k(\bt')}\,
\bbu(
\textrm{$\bt'$ satisfies the conditions (1) and (2) above})
\\
&=
\frac{1}{A_n}\, 
\frac{\prod_j W(t'_j)}{k(\bt')}\,
\delta\Big(\sum_j v_j-n\Big)
\,
\delta\Big(\sum_j \ell_j-k+1\Big)
\ef,
\end{split}
\ee
with complexity $T(n)$, we are able to exactly sample our
combinatorial structures with complexity $T(n) + \Theta(n)$.

In the remaining of this paper we will concentrate on the problem of
exactly sample lists $\bt$ from the measure $\mu_n^{\rm (cyc)}(\bt)$
above. At this point, the pertinence of our toy-model problem of
Section \ref{sec.simpleex}, dealing with the simplest possible family
of lattice bridges which is not Dyck bridges (for which the problem
reduces to random shhuffles), should be apparent.

\section{Further analytic preliminaries}
\label{sec.anaprel}

\noindent
Having renamed $A(z)$ our generating function $Y^{(0)}(z)$, our system
(\ref{eq.2664576743}) reads
\be
\left\{
\begin{array}{l}
Y^{(0)}(z) \equiv A(z)
= \Phi^{(0)}(z,A(z),Y^{(1)}(z),\ldots,Y^{(m)}(z))
\\
%
Y^{(1)}(z)
= \Phi^{(1)}(z,A(z),Y^{(1)}(z),\ldots,Y^{(m)}(z))
\\
\setlength{\unitlength}{8pt}
\begin{picture}(0,1)
\put(1.8,1.3){$\scriptstyle{\vdots}$}
\end{picture}%
\rule{0pt}{22pt}%
Y^{(m)}(z)
= \Phi^{(m)}(z,A(z),Y^{(1)}(z),\ldots,Y^{(m)}(z))
\end{array}
\right.
\ee
Similarly the generating function for the subtrees
\be
A(z,u) = 
\sum_{t \in \cT^{\rm (sub)}} W(t)\, z^{v(t)} u^{\ell(t)}
=
\sum_{v, \ell \geq 0}
z^v u^{\ell} 
\sum_{t \in \cT_{v,\ell}} W(t)
\ee
is given by
\be
A(z,u)
= \Phi^{(0)}(z,u,Y^{(1)}(z,u),\ldots,Y^{(m)}(z,u))
\ee
where the $Y^{(\alpha)}(z,u)$'s satisfy the reduced system 
\be
\left\{
\begin{array}{l}
Y^{(1)}(z,u)
= \Phi^{(1)}(z,u,Y^{(1)}(z,u),\ldots,Y^{(m)}(z,u))
\\
\setlength{\unitlength}{8pt}
\begin{picture}(0,1)
\put(1.8,1.3){$\scriptstyle{\vdots}$}
\end{picture}%
\rule{0pt}{22pt}%
Y^{(m)}(z,u)
= \Phi^{(m)}(z,u,Y^{(1)}(z,u),\ldots,Y^{(m)}(z,u))
\end{array}
\right.
\ee
As a check, let us verify that 
the measure in (\ref{eq.37645734c}) is indeed normalised.
At this aim we should have, for all $n \geq 0$,
\be
\label{eq.37645734d}
\begin{split}
1 = \sum_{\bt} \mu_n^{\rm (cyc)}(\bt) &= 
\frac{1}{A_n}\, 
[z^n u^{-1}]
\sum_k
\frac{1}{k}
\Big(
\sum_{t} W(t) z^{v(t)} u^{\ell(t)-1}
\Big)^k
\\
&=
-\frac{1}{A_n}\, 
[z^n u^{-1}]
\ln
\Big(
1-\frac{A(z,u)}{u}
\Big)
\\
&=
-\frac{1}{A_n}\, 
[z^n]
\oint
\frac{\dxns{u}}{2 \pi i}
\ln
\Big(
1-\frac{A(z,u)}{u}
\Big)
\ef,
\end{split}
\ee
where the integral is on a contour encircling the origin, and no other
pole of the integrand. Or, in other words, taking the generating
function for an arbitrary value of $z$ in the interval $[0,\rho[$,
\be
\label{eq.37645734d}
\begin{split}
A(z) 
&= 
-
\oint
\frac{\dxns{u}}{2 \pi i}
\ln
\Big(
1-\frac{A(z,u)}{u}
\Big)
\ef.
\end{split}
\ee
However, the logarithm has a cut discontinuity on an interval starting
from the origin, along the positive real axis, and this cut ends at
the point $u(z)$ such that $A(z,u(z))=u(z)$. If we deform the countour
of integration as to encircle this cut discontinuity, the integral is
exactly (minus) the length of the interval (factors $2 \pi i$ simplify
exactly as in the classical proof of the Cauchy Theorem), so we get
that our check is equivalent to the relation $A(z)=u(z)$ for all 
$z \in [0,\rho[$, whenever $A(z,u(z))=u(z)$. Indeed, on the manifold
$(z,u) \in \cM \subset \bC^2$ where the latter holds, any solution of
the bivariate reduced system is also a solution of the original
univariate system, i.e.\ $Y^{(\alpha)}(z,u(z))=Y^{(\alpha)}(z)$ for
all $0\leq \alpha \leq m$, and in particular for $m=0$, so that our
claim is verified.

There is a two-parameter family of natural measures on the set of
possible subtrees $t$, where the two parameters are the Lagrange
multipliers for the number of $A$- and $z$-leaves:
\be
\mu(t;z,u):=\left(\frac{A(z,u)}{u}\right)^{-1}
z^{v(t)} u^{\ell(t)-1}
\ef.
\ee
One point $(\zeta,\theta)\in (\bR^+)^2$, that we shall call the
\emph{critical point}, has a crucial role in the following. These
two values are defined by
\begin{align}
\label{eq.condZT}
\frac{A(\zeta,\theta)}{\theta}
&=
1
\ef;
&
\deenne{\theta}{{}}
\frac{A(\zeta,\theta)}{\theta}
&=
0
\ef.
\end{align}
The first condition corresponds to $(\zeta,\theta) \in \cM$, i.e., to
the fact that the weight $\theta$ associated to the $A$-leaves
coincides with the generating function of trees with root-index
$A$. The second condition is a rewriting of the equation
$\deenne{\theta}{{}} \ln \frac{A(\zeta,\theta)}{\theta}=0$
in light of $(\zeta,\theta) \in \cM$, and this equation describes the
fact that our generalised \Luk\ paths have average zero drift, as to
be expected for optimising the efficiency of the Boltzmann sampling of
paths from $(0,0)$ to $(n,-1)$, when $n \gg 1$.  If we had to perform
a Boltzmann sampling of our generalised \Luk\ bridges, the optimality
strategy suggested by singularity analysis would tell that we have to
sample the subtrees $t_j$ of the walk exactly with the measure
$\mu(t;\zeta,\theta)$. However, this first algorithm would have a
rather poor efficiency: not only we have a factor $\sim 1/\sqrt{n}$
for the probability of reaching the point $(n-1)$, but also we have a
crude acceptance rate of order $\sim 1/n$ on the typical cases,
because our acceptance rate must reproduce the factor $1/k(\bt)$ in
the cyclic lemma up to a multiplicative constant, but this factor, in
general, may be as large as 1 (in the exponentially rare case in which
the root is the only node of the tree with label $A$). Note that we
did not have this problem in our toy model of Section
\ref{sec.simpleex}, as in that case the possible steps of the walk
have bounded length.  As we will see, we will need to come up with two
different ideas, for dealing with these two different problems.
\label{pag.twoprobl}

There is a second remarkable characterisation of the critical point
$(\zeta,\theta)$, which is the following:
\begin{lemma}
The values $(\zeta,\theta)$ satisfying (\ref{eq.condZT}) are the
components $(\rho,\tau_0)$ of the solution $(\rho,\vec{\tau})$ of the
characteristic system (\ref{eq.charsys}). Furthermore, the vector 
$(1,\deenne{\theta}{{}}Y^{(1)}(\zeta,\theta)),\ldots,\deenne{\theta}{{}}Y^{(m)}(\zeta,\theta))$
is the Frobenius vector of the matrix $K$.
\end{lemma} 

\noindent
{\it Proof.}
Indeed, the first equation of (\ref{eq.condZT}) coincides with the
component $0$ of the first equation of (\ref{eq.charsys}). 
For the second equation of (\ref{eq.condZT}), 
the explicit calculation of $\deenne{\theta}{{}}A(\zeta,\theta)$ gives
\be
\begin{split}
\deenne{\theta}{{}}A(\zeta,\theta)
&=
\dienne{\theta}{{}}\Phi^{(0)}(\zeta,\theta,
 Y^{(1)}(\zeta,\theta),\ldots,Y^{(m)}(\zeta,\theta))
\\
&
=
\deenne{\theta}{{}}\Phi^{(0)}(\zeta,\theta,
 Y^{(1)}(\zeta,\theta),\ldots,Y^{(m)}(\zeta,\theta))
\\
& \quad
+
\sum_{\beta=1}^m
\deenne{{Y^{(\beta)}}}{{}}\Phi^{(0)}(\zeta,\theta,
 Y^{(1)}(\zeta,\theta),\ldots,Y^{(m)}(\zeta,\theta))
\deenne{\theta}{{}}Y^{(\beta)}(\zeta,\theta)
\\
& \qquad
=
K_{00}
+
\sum_{\beta=1}^m
K_{0\beta}
{\deenne{\theta}{{}}}
Y^{(\beta)}(\zeta,\theta)
\ef.
\end{split}
\ee
Similarly, for the other components we have
\be
\begin{split}
\deenne{\theta}{{}}
Y^{(\alpha)}(\zeta,\theta)
&=
K_{\alpha 0}
+
\sum_{\beta=1}^m
K_{\alpha \beta}
\deenne{\theta}{{}}Y^{(\beta)}(\zeta,\theta)
\ef.
\end{split}
\ee
That is, introducing the shortcuts
\begin{align}
\vec{v}
&=
\Big(
\deenne{\theta}{{}}A(\zeta,\theta),\deenne{\theta}{{}}Y^{(1)}(\zeta,\theta),\ldots,\deenne{\theta}{{}}Y^{(m)}(\zeta,\theta)
\Big)
\ef;
\\
\vec{v}'
&=
\Big(
1,\deenne{\theta}{{}}Y^{(1)}(\zeta,\theta),\ldots,\deenne{\theta}{{}}Y^{(m)}(\zeta,\theta)
\Big)
\ef;
\end{align}
we have
\be
K \vec{v}'=\vec{v}
\ef.
\ee
If both equations (\ref{eq.condZT}) hold, we have
\be
0
=
\deenne{\theta}{{}}
\frac{A(\zeta,\theta)}{\theta}
=
\frac{\theta \deenne{\theta}{{}} 
A(\zeta,\theta) - A(\zeta,\theta)
}{\theta^2}
=
\frac{1}{\theta} 
\left(\deenne{\theta}{{}} A(\zeta,\theta) - 1\right)
\ef,
\ee
that is, $\vec{v}=\vec{v'}$, which thus is an eigenvector of $K$, with
eigenvalue 1. Now, the fact that the coefficients $\phi$ defining the
system are real-positive implies that $K$ has a Frobenius eigenvector,
and also that the entries of $\vec{v}$ are all positive, so 
$\vec{v}$ must be the Frobenius eigenvector of $K$, and
$(\zeta,\theta,Y^{(1)}(\zeta,\theta),\ldots,Y^{(m)}(\zeta,\theta))$
coincides with the solution $(\rho,\tau_0,\tau_1,\ldots\tau_m)$ of the
characteristic system.
\hfill $\square$

\medskip
\noindent
One consequence of this lemma is that, in lack of a more direct
strategy specific for the structure at hand, we can always determine
$(\zeta,\theta)$ by using general-purpose libraries already developed
for the ordinary Boltzmann Method, as explained in
\cite{PivSalSor2012}, or, when we will deal with multi-parametric
Boltzmann sampling, in \cite{SergeyMultipar} (see also 
{\tt https://paganini.readthedocs.io/en/latest/}).

A further analytic ingredient that we will need in the following is:
\begin{lemma}
\label{lem.ballZT}
There exists a radius $\delta>0$ such that the measure
$\mu(t;z,u)$ is associated to a subcritical GW process, for all $(z,u)
\in \cB_{(\zeta,\theta), \delta}$.
\end{lemma}

\noindent
{\it Proof.}  The rewriting process associated to the original system
$\vec{Y}(z)=\vec{\Phi}(z,\vec{Y}(z))$ gives a critical GW process when
$(z,\vec{Y})=(\rho,\vec{\tau})$ satisfies the equations
(\ref{eq.charsys}), that is, as implied by the previous lemma, when
$(z,\vec{Y})=(\zeta,\theta,Y^{(1)}(\zeta,\theta),\ldots,Y^{(m)}(\zeta,\theta))$.
Now, it is intuitively clear that, if in a critical rewriting system
we change one rewriting rule with a stopping rule, the new system is
sub-critical. A formal proof is given by Lemma \ref{lem.frobAB}, in
Appendix \ref{app.frobfact}, specialised to $x=1$ (for the original
system) and $x=0$ (for the reduced system). This implies the statement
for the point $(\zeta,\theta)$. The continuity of the Frobenius
eigenvalue as a function of $z$ and $u$ implies the existence of
a non-empty neighbourhood centered in $(\zeta,\theta)$, and completes
the proof.
\hfill $\square$

\section{The algorithm}
\label{sec.algo}

\noindent
Let us define the support $U$ of the combinatorial specification as
the set of vectors 
$\vec{u}=(v,\ell-1) \in \{0,1,\ldots\} \times \{-1,0,1,\ldots\}$
such that $\cT_{v,\ell} \neq \varnothing$. Note that $(0,0)$ and
$(0,-1)$ cannot be elements of $U$ (because every tree must have some
leaves, and we have required the linear part of $\vec{\Phi}$ to be
nilpotent, so we cannot have only one leaf, and of $A$ type).
Let us also call $v_0 \geq 1$ the
smallest integer such that $(v_0,-1) \in U$, i.e.\ such that
$\cT_{v_0,0} \neq \varnothing$. This integer must exist, otherwise
there would exist no generalised \Luk\ excursions associated to our
combinatorial class.  
Call
\be
T_0 = \sum_{t \in \cT_{v_0,0}} W(t)
\ef,
\ee
the generating function of our objects in $\cT_{v_0,0}$ (which is a
finite polynomial in the $phi$'s), and call $\mu_0(t)$ the natural
measure on $\cT_{v_0,0}$
\be
\label{eq.mu0}
\mu_0(t)=\frac{W(t)}{T_0}
\ee
(note that it is also the restriction of $\mu(t;\zeta,\theta)$ to
$\cT_{v_0,0}$, although, in fact, it does not depend on $\zeta$ and
$\theta$).

Let $\mu_{\neq}(t)$ be some measure, still to be determined, with
support on $\cT^{\rm (sub)} \smallsetminus \cT_{v_0,0}$.

Let us call the \emph{landing diagonal} $D_n$ as the set of points in 
$\bN \times \bZ$ on the line with slope $-1/v_0$ passing through
$(n,-1)$, and with ordinate above $-1$:
\be
D_n = \{ (i,j) \;|\; j \geq -1, \; i+v_0 j = n-v_0 \}
\ef.
\ee
Let us call the 
\emph{final region} $F_n$ as the set of points in 
$\bN \times \bZ$ beyond the landing diagonal:
\be
F_n = \{ (i,j) \;|\; i+v_0 j \geq n-v_0 \}
\ef.
\ee
Our algorithm for random (generalised \Luk) bridges induced by an
irreducible context-free structure has the same idea outlined in our
toy model of Section \ref{sec.simpleex}, that is, for a suitable
measure $\mu_{\neq}(t)$ and a suitable acceptance-rate function
$r_n(\{\vec{u}_j\}_{1 \leq j \leq k(\bt)})$,

\begin{algorithm}[H]
\Begin{
\Repeat{$\vec{u} \in D_n$ and $\mathrm{rand}_{[0,1]}<r_n(\{\vec{u}_j\})$}
{
$\vec{u}=(0,0)$\;
$k=0$\;
$\bt=(\,)$\;
\While{$\vec{u} \not\in F_n$}
{
 $k \to k+1$\;
 $t_k \longleftarrow \mu_{\infty}$\;
 $\vec{u}_k=(v(t_k),\ell(t_k)-1)$\;
 $\vec{u} \to \vec{u} + \vec{u}_k$\;
 Append $t_k$ to $\bt$\;
}
}
$m=\vec{u}_2+1$\;
\For{$i=1\ldots m$}
{ 
 $t_{k+i} \longleftarrow \mu_{0}$\;
 Append $t_{k+i}$ to $\bt$\;
}
$\sigma \longleftarrow \textrm{BBHL}(k,m)$\;
$\bt = 
\sigma \circ \bt$\;
\Return{$\bt$}
}
\caption{Accelerated Boltzmann Method, for generalised \Luk\ bridges.\label{algo.my}}
\end{algorithm}
\noindent
The algorithm for excursions, and thus for (coloured) trees and for
irreducible context-free structures, is identical to the previous one, in
which the last instruction
$\bt = \sigma \circ \bt$
is replaced by
$\bt = \cyc \circ \sigma \circ \bt$.

At this point we have to guess which measure $\mu_{\neq}(t)$ is more
adapted at this aim.  Let us call $\mu_{\neq}(t;z,u)$ the restriction
of $\mu(t;z,u)$ to $\cT \smallsetminus \cT_{v_0,0}$.  If we believe in
the Boltzmann Method paradigm, we should be induced to choose, as
measure $\mu_{\neq}(t)$, the critical measure
$\mu_{\neq}(t;\zeta,\theta)$.  However, in the past, experience has
shown that small modifications of the critical parameters, scaling
algebraically with $n$, may improve the performance of the algorithm
(see the papers cited in the introduction, or also,
e.g.\ \cite[Sec.\;4.4]{PivoPP}), and, as we will see later on, this is
the case also here.

Let us then introduce a larger family of measures.  At fixed $n$, let
$X=\{ (\pi_a,z_a,u_a) \}_{1 \leq a \leq q}$ be a finite set of
parameters, where $\pi_a$ is a discrete probability distribution
($\pi_a>0$ for all $a$, and $\pi_1+\cdots+\pi_q=1$), and $(z_a,u_a)$
are points in the ball introduced in Lemma \ref{lem.ballZT}. We will
define the measure $\mu_{\neq}(t)$ associated to $X$ as 
\be
\label{mumix}
\mu_{\neq}(t) = \sum_{a=1}^q \pi_a \mu_{\neq}(t;z_a,u_a) 
\ef.
\ee
Sampling from these measures is immediate: you should just first
sample $a$, with distribution $\pi$, and then sample from
$\mu(t;z_a,u_a)$, using the stochastic rewriting system, possibly
repeating the procedure (with the same $a$) if the subtree $t$ is in
$\cT_{v_0,0}$ (both this event and its complementary happen with a
probability $\Theta(1)$, so the average number of runs is also 
$\Theta(1)$). 
More explicitly, the measure above reads 
\be
\mu_{\neq}(t)
=
\sum_{a=1}^q
\pi_a
\frac{u_a}{A(z_a,u_a) - T_0\, z_a^{v_0}}
W(t)
z_a^{v(t)}
u_a^{\ell(t)-1}
=:
\sum_{a=1}^q
\tilde{\pi}_a
\frac{\theta}{A(\zeta,\theta) - T_0\, \zeta^{v_0}}
W(t)
z_a^{v(t)}
u_a^{\ell(t)-1}
\ef,
\ee
where the parameters $\tilde{\pi}_a$ have been introduced for later
convenience, but are in general not normalised.

We can reformulate
a set of parameters $X$ 
in the form $z_a = \zeta \exp(\xi_a)$, $u_a = \theta \exp(\eta_a)$..
Call $\bar{\xi}= \sum_a \tilde{\pi}_a \xi_a/\sum_a \tilde{\pi}_a$
and $\bar{\eta}=\sum_a \tilde{\pi}_a \eta_a/\sum_a \tilde{\pi}_a$. 
In this reparametrisation we have
\be
\label{eq.6486534674}
\mu_{\neq}(t)
=
\frac{\theta}{A(\zeta,\theta) - T_0\, \zeta^{v_0}}
W(t)
\zeta^{v(t)}
\theta^{\ell(t)-1}
\sum_{a=1}^q
\tilde{\pi}_a
\exp(\xi_a v(t)+\eta_a (\ell(t)-1))
\ef.
\ee
We say that a set $X$ is \emph{balanced} if $\bar{\xi}= \bar{\eta}=0$.
Of course, the critical sampler $X=\{(1,\zeta,\theta)\}$ is also
balanced. 

In the second part of the algorithm we must repeatedly sample from
$\mu_0$. This is algorithmically is immediate, as in fact the set
$\cT_{v_0,0}$ contains finitely many elements, which can be found in
constant time in a preprocessing procedure, together with their
weights $W(t)$.

So, all that is left to do is repeating our analysis of the
acceptance-rate function, along the same lines presented in Section
\ref{sec.simpleex}, but in this more complicated setting.  We do this
analysis in the next sections.

\section{Analysis of the acceptance ratio}
\label{sec.anAR}

\noindent
The conditions on the acceptance-rate function read in this
case:
\begin{enumerate}
\item the algorithm samples according to the desired measure;
\item $r_n(\{\vec{u}_j\}) \in [0,1]$ for all $n, k\geq 1$ and
  lists $\{\vec{u}_j\}_{1 \leq j \leq k}$ with $\sum_j \vec{u}_j \in D_n$;
\item the function $r_n(\{\vec{u}_j\})$ can be calculated
  efficiently, i.e., in average time at most $\cO(n)$;
\item $\eee(r_n(\{\vec{u}_j\}))=\Theta(1)$.
\end{enumerate}
\noindent
Again, we start by analysing the probability of getting any given
output string $\bt$. Recall that we have called $m$ the integer
corresponding to the ordinate of the point on the landing diagonal
that we reached, plus 1. In other words, if we label the points of the
landing diagonal with the non-negative integers, we have landed on the
$m$-th point. We also called $k$ the number of subtrees generated in
the first part of the algorithm, thus the complete list $\bt$ obtained
at the end of the algorithm has $k+m$ items, and, for $n>v_0$, has
$1 \leq k\leq n-v_0$ and $0 \leq m \leq \lfloor n/v_0 \rfloor-1$. 
Indeed, the points $\vec{u}$ of $\bN \times \bZ$ are partitioned into diagonals
according to the integer 
$\tilde{u}(\vec{u}) := u_1 + v_0 u_2 = \vec{u} \cdot (1,v_0)$. 
By definition of $v_0$, for all the points $\vec{u}$ in the support
$U$, this combination is a positive integer, except for $(v_0,-1)$,
for which it is zero. The origin is on the diagonal with parameter
zero, and the landing diagonal has parameter $n':=n-v_0$. The case with
largest possible values of $k$ and $m$ is when only the steps $(0,1)$
and $(v_0,-1)$ are used, in the first and second part of the
algorithm, respectively, which occurs, for example, in the case of
binary trees counted by number of leaves, $A=z+A^2$ (which, from the
point of view of our algorithm, is trivial, as the core of the
operations is fully contained in the BBHL shuffling procedure).
We anticipate that we will only consider functions
$r_n(\{\vec{u}_j\})$ that only depend on the parameters 
$\tilde{u}_j=\tilde{u}(\vec{u}_j)$.

Introduce the shortcuts
\begin{align}
A_0 
&= \frac{T_0\, \zeta^{v_0}}{\theta}
\ef;
&
A_{\neq}
&= 
\frac{A(\zeta,\theta)}{\theta}-A_0
\ef.
\end{align}
Call
\be
\bar{v}:= \zeta \deenne{\zeta}{{}} \frac{A(\zeta,\theta)}{\theta}
\ef.
\ee
This quantity is the average horizontal length of the steps of our
\Luk\ walks, if these are sampled with the critical measure
$\mu(t;\zeta,\theta)$.
That is, for critical steps, 
\be
\eee(\vec{u})_{\mu(t;\zeta,\theta)} = (\bar{v},0)
\ef.
\ee
As the GW process associated to the subtrees
is subcritical, this quantity is of order $1$.

It is also easy to calculate the analogous quantity for
$\mu_{\neq}(t;\zeta,\theta)$, as
\be
(\bar{v},0)
=
\eee(\vec{u})_{\mu(t;\zeta,\theta)} 
= 
A_{\neq}
\eee(\vec{u})_{\mu_{\neq}(t;\zeta,\theta)}
+
A_0
\eee(\vec{u})_{\mu_0(t)}
=
A_{\neq}
\eee(\vec{u})_{\mu_{\neq}(t;\zeta,\theta)}
+
A_0
(v_0,-1)
\ef,
\ee
that is
\be
\eee(\vec{u})_{\mu_{\neq}(t;\zeta,\theta)}
=
\frac{(\bar{v}-A_0 v_0,A_0)}{1-A_0}
\ef.
\ee
This has two important consequences. First, the initial part of the
algorithm makes walks with a positive drift. Second, the average
amount of diagonals by which these walks do jump at each step is
$\frac{(\bar{v}-A_0 v_0,A_0)}{1-A_0} \cdot (1,v_0) =
\bar{v}/A_{\neq}=\Theta(1)$.  This is an aspect of the fact that the
probability of reaching the landing diagonal is of order $1$. More
precisely, if $p(\tilde{u})$ is the probability distribution of
$\tilde{u}$ for steps in a walk, and
$\hat{p}(x)=\sum_{\tilde{u}}p(\tilde{u}) x^{\tilde{u}}$ is the
corresponding generating function, the asymptotic probability of
reaching the landing diagonal is 
$\left.(\deenne{x}{{}} \hat{p}(x))^{-1}\right|_{x=1}$,
and the corrections are exponentially
small in the length of the walk.\footnote{More
  generally, in a Markov chain in which we jump by $a$ with
  probability $p_a$ (with $a \in \bN^+$), the asymptotic probability
  $q$ that a site $n$ is not occupied is given by the sum over 
  $k,j \geq 0$ of the probability that we have a jump of length
  $k+j+2$, from $n-k-1$ (which must be occupied), to $n+j+1$:
\be
q=\sum_{j,k \geq 0} (1-q)\, p_{j+k+2}
=\sum_{\ell \geq 2} (1-q)\, (\ell-1) p_{\ell}
= (1-q) (\hat{p}'(1)-1)
\ee
that is, the probability that a site is occupied is
$1-q=1/\hat{p}'(1)$.
}
So, in our case,
\be
\label{eq.ContoPDn}
\ppp(\textrm{reaching $D_n$})
=
\left. \left(
\deenne{x}{{}}
\frac{
x^{-v_0}
A(\zeta x,\theta x^{v_0})
-
\theta 
A_0
}{
A(\zeta,\theta)-\theta A_0
}
\right)^{-1}\right|_{x=1}
+
\cO(\exp(-\alpha n))
=
\Theta(1)
\ef.
\ee
This is an important check, as 
$\ppp(\textrm{reaching $D_n$})=\Theta(1)$ is a necessary condition for
an algorithm of our form to reach complexity optimal up to a
multiplicative constant. The inverse of the quantity evaluated above
enters as a factor in this complexity multiplicative constant.

Now, let $\bt$ be a list that, in the first part of our algorithm, has
reached the landing diagonal. The probability of getting this, with
parameters $k$ and $m$ as defined above, before the shuffling, is
\be
p(\bt)
=
\Big(
\prod_{j=1}^{k}
\mu_{\neq}(t_j)
\Big)
\Big(
\prod_{j=1}^{m}
\mu_{0}(t_{k+j})
\Big)
r_n(\{\tilde{u}_j\})
\ef.
\ee
After the shuffling, we get a further factor $1/\binom{k+m}{k}$.
The resulting probability must be proportional to
the measure $\mu_n^{\rm (cyc)}(\bt)$ in
(\ref{eq.37645734c}), that is
\be
\label{eq.53787536545}
\begin{split}
\Big(
\prod_{j=1}^{k}
\mu_{\neq}(t_j)
\Big)
\Big(
\prod_{j=1}^{m}
\mu_{0}(t_{k+j})
\Big)
r_n(\{\tilde{u}_j\})
\propto
\binom{k+m}{k}
\frac{\prod_j W(t_j)}{k+m}\,
\ef.
\end{split}
\ee
Note that we can rewrite (\ref{eq.6486534674}) as
\begin{align}
\label{eq.6486534674b}
\mu_{\neq}(t)
&=
\frac{W(t)
\zeta^{v(t)}
\theta^{\ell(t)-1}
}{A_{\neq}}
F(v(t),\ell(t))
\ef;
&
F(v,\ell)
&=
\sum_{a=1}^q
\tilde{\pi}_a
\exp(\xi_a v+\eta_a (\ell-1))
\ef.
\end{align}
More simply, (\ref{eq.mu0}) reads
\be
\label{eq.mu0b}
\mu_0(t)=\frac{W(t) \zeta^{v_0} \theta^{-1}}{A_0}
\ef,
\ee
thus, substituting in (\ref{eq.53787536545}),
\be
\label{eq.53787536545b}
\begin{split}
\zeta^{n} \theta^{-1}
\Big(
\prod_{j=1}^{k+m} 
W(t_j)
\Big)
\Big(
\prod_{j=1}^{k}
F(v_j,\ell_j)
\Big)
r_n(\{\tilde{u}_j\})
\propto
\frac{
(k+m-1)!}{k!\, m!}
A_{\neq}^k A_0^m
\prod_{j=1}^{k+m} 
W(t_j)
\ef.
\end{split}
\ee
That is, there must exist a function $K_n$ such that
\be
\label{eq.53787536545b}
\begin{split}
r_n(\{\tilde{u}_j\})
=
K_n
\frac{
(k+m-1)!}{k!\, m!}
A_{\neq}^k A_0^m
\Big(
\prod_{j=1}^{k}
F(v_j,\ell_j)^{-1}
\Big)
\ef,
\end{split}
\ee
which can be calculated efficiently, and such that $r$ is valued in
$[0,1]$, and its average is $\Theta(1)$.

In particular, we have required that $r_n$ is a function of the
$\tilde{u}_j$'s only, that is, we must have that $F$ is a function of
$\tilde{u}$ only. This happens if and only if $\eta_a=v_0 \xi_a$ for
all $a$, and in this case we have
\begin{align}
\label{eq.6486534674bx}
F(\tilde{u})
&=
\sum_{a=1}^q
\tilde{\pi}_a
\exp(\xi_a \tilde{u})
\ef.
\end{align}
We have anticipated that ``critical sampling'', that is, performing
the sampling in the first part of the algorithm with the critical
measure $\mu_{\neq}(t)=\mu_{\neq}(t;\zeta,\theta)$, does not lead to
an optimal complexity in general.  Nonetheless, it is instructive to
first analyse this simpler case, and understand the reason why we need
the more complicated combinations in equation (\ref{mumix}).

In the critical case, the function $F(\tilde{u})$ is just $1$, and we have
\be
\label{eq.53787536545c}
\begin{split}
r_n(\{\tilde{u}_j\})
=
K_n
\frac{
(k+m-1)!}{k!\, m!}
A_{\neq}^k A_0^m
\ef.
\end{split}
\ee
Recall that, from the definition of $A_{\neq}$ and $A_0$, and the
equation (\ref{eq.condZT}), we have $A_{\neq}+A_0=1$, so that the
combination $\frac{(k+m)!}{k!\, m!}  A_{\neq}^k A_0^m$ is a binomial
distribution. At fixed $k+m$, it is maximised at 
$k/(k+m) \simeq A_{\neq}$. The further factor coming from the cyclic
lemma gives only a small perturbation.
Analytic Combinatorics predicts that, for lists $\bt$ which are
typical in the seeked measure $\mu_n^{\rm (cyc)}(\bt)$, we have
$(k+m) \bar{v} \simeq n$. Thus, Stirling approximation gives
$r_n(\{\tilde{u}_j\}) \sim K_n C/n^{\frac{3}{2}}$ on typical lists, with
$C$ a finite constant that can be easily evaluated.

However, this result is not satisfactory. Indeed, unless the reduced
system is specially simple, there exist trees $T$ in $\cT$ with root
of label $A$, no other node with this label, and arbitrarily large
number of $z$-leaves.  In this case, we can have lists $\bt(T)$ of
subtree decompositions with values of $(k,m)$ as small as $(1,0)$, and
in this case we have
\be
r_n(\{n-v_0\}) = 
K_n A_{\neq}
\ef.
\ee
As we must have $\max_{\{\tilde{u}_j\}} r_n(\{\tilde{u}_j\}) \leq 1$, we
have $K_n \leq 1/A_{\neq}$, and the acceptance rate on typical
lists is only $\cO(n^{-\frac{3}{2}})$. This is the `second problem' of
the two mentioned on page \pageref{pag.twoprobl}.

This suggests to explore the more general class of measures in
equation (\ref{mumix}), where the function $F(\tilde{u})$ may enter the
expression (\ref{eq.53787536545b}) in a way that solves this
problem. In the following we shall consider only balanced sets
$X$.\footnote{Indeed, it can be seen that non-balanced sets $X$ lead
  to weights of the form $\exp(\bar{\xi}(n-v_0 m)+\bar{\eta}(m-1)+c
  k)$ times the one associated to the balanced set, in which the
  $\xi_a$'s and $\eta_a$'s are translated by their weighted averages,
  and the $\pi_a$'s are modified as to produce the same set
  $\tilde{\pi}_a$ up to an overall factor $e^{c}$, and factors of this
  form are not the best choice for correcting the $\sim
  n^{-\frac{3}{2}}$ term in the acceptance rate.}
The conceptually simplest choice is to take a set $X$ composed of only
two points, $a \in \{+,-\}$, with $\tilde{\pi}_+=\tilde{\pi}_-$ and
$\xi_+ = -\xi_-$, $\eta_+ = -\eta_-$. We anticipate that the
appropriate scaling is of the form $\xi_{\pm}, \eta_{\pm} =\Theta(n^{-\frac{1}{2}})$,
and, for definiteness, we choose once and for all
\begin{align}
\xi_{\pm} &= 
\pm \frac{\eps}{\sqrt{n-v_0}}
\ef;
&
\eta_{\pm} &= 
\pm \frac{\eps}{\sqrt{n-v_0}}
v_0
\ef;
\end{align}
with $\eps$ a constant to be determined later on. (In fact, optimising
$\eps$ is not crucial, as it only affects the overall constant in the
complexity, not the scaling exponent.)  Of course, using $n$ and
$n'=n-v_0$ is essentially equivalent for large $n$, and our choice is
driven by the fact that the combination $n'$, the index of the landing
diagonal, appears repeatedly in the analysis.

This scaling in $n$ implies, by continuity, that the considerations
above (e.g.\ on the drift of the first part of the walk), leading to
the fact that we have a positive probability of reaching the landing
diagonal, remain valid in this more general framework, by continuity,
for $n$ large enough. Also, for $n$ large enough, the points
$(z_{\pm},u_{\pm})$ are contained in the disk of Lemma \ref{lem.ballZT}.

Consider the function
\be
f(x) =
\ln 
\Bigg(
\frac{
e^{-x v_0}
A(\zeta e^{x},\theta e^{x v_0})
-A_0
}
{
A(\zeta,\theta)
-
A_0
}
\Bigg)
\ef.
\ee
We have 
$\tilde{\pi}_{\pm} = \pi_{\pm} \exp(\pm f(\frac{\eps}{\sqrt{n'}}))$. As
we have required $\pi_+ + \pi_- = 1$ and
$\tilde{\pi}_+ = \tilde{\pi}_-$, we have
\be
\tilde{\pi}_+ + \tilde{\pi}_-
=
\frac{1}{\cosh(f(\frac{\eps}{\sqrt{n'}}))}
=
1-\Theta(n^{-1})
\ef.
\ee
As a result, the factors $F(\tilde{u})$ appearing in
(\ref{eq.53787536545b}) take the simple form
\begin{align}
F(\tilde{u})
&=
\frac{\cosh(\frac{\eps}{\sqrt{n'}}
\tilde{u})}
{\cosh(f(\frac{\eps}{\sqrt{n'}}))}
\ef.
\end{align}
From the convexity of the hyperbolic cosine (that is, 
$\min_{p} \cosh(px) \cosh((1-p)x)$ is attained at $p=\frac{1}{2}$)
we have that, for any value of $k$, the maximum of the combination
$\prod_{j=1}^{k} F(\tilde{u}_j)^{-1} $ appearing in
(\ref{eq.53787536545b}) is bounded from above by the $k$-th power of
the analytic continuation of the function $F$, valued at 
$k^{-1} \tilde{u}(n,-1)=n'/k$, that is
\be
\label{eq.53787536545d}
\begin{split}
r_n(\{\tilde{u}_j\})
\leq
K_n
\frac{
(k+m-1)!}{k!\, m!}
A_{\neq}^k A_0^m
\left(
\frac
{\cosh(f(\eps \cdot \frac{1}{\sqrt{n'}}))}
{\cosh(
\eps \cdot
\frac{\sqrt{n'}}{k}
)}
\right)^{k}
\ef.
\end{split}
\ee
As $k$ is bounded by $n'$, 
for all lists $\bt$
the factor 
$\cosh(f(\eps \cdot \frac{1}{\sqrt{n'}}))^k$
is in a range $[1,C]$, for some $C=\Theta(1)$ that can be evaluated easily, 
and we can omit this factor in the analysis of an upper bound (up to
include a factor $C$ at the end of the calculation), i.e.\ search a
value $K_n$ that satisfies the bound
\be
\label{eq.53787536545d2}
\begin{split}
r_n(\{\tilde{u}_j\})
\leq
K_n
C
\frac{
(k+m-1)!}{k!\, m!}
A_{\neq}^k A_0^m
\cosh
\Big(
\frac{\eps \sqrt{n'}}{k}
\Big)^{-k}
\ef.
\end{split}
\ee
Now, if $k \sim n^{\gamma}$
we have
\be
\left(
\cosh(\eps \sqrt{n'}/k)
\right)^{-k}
\sim
\left\{
\begin{array}{ll}
\rule{0pt}{2pt}%
\exp(-\Theta(\sqrt{n})) & \gamma \leq \frac{1}{2} \\
\rule{0pt}{12pt}%
\exp(-\Theta(n^{1-\gamma})) & \frac{1}{2} < \gamma <1
\\
\rule{0pt}{12pt}%
\exp(-\eps^2/(2 \kappa))
& 
\gamma=1, \; k \simeq \kappa n
\end{array}
\right.
\ee
so, as for lists with $k \sim n$ we have a bound of order
$n^{-\frac{3}{2}}$, while 
for lists with $k \sim n^{\gamma} = o(n)$ we have a bound not larger
than $\exp(-\Theta(n^{1-\gamma}))$, we can perform our analysis under
the assumption that $k=\Theta(n')$.
At $k$ fixed and large, the maximum in $m$ 
is achieved for 
$m=\lfloor(k A_0-1)/A_{\neq}\rfloor$, we have that also $m=\Theta(n')$. 

So we have proven that the worst case of the bound is achieved for 
$k=\kappa n'$, $m=\mu n'$, and $\kappa, \mu = \Theta(1)$.
Furthermore, the saddle-point approach predicts
$\kappa/\mu = A_{\neq}/A_0 + o(1)$, and the bound above reads
\be
\label{eq.53787536545d3}
\begin{split}
r_n(\{\tilde{u}_j\})
\lesssim
K_n
C
\frac{1}{n^{\frac{3}{2}}}
\frac{A_{\neq}}{\sqrt{2 \pi A_0}}
\frac{1}{\kappa^{\frac{3}{2}}}
\exp\left(-\frac{\eps^2}{2 \kappa}\right)
\ef.
\end{split}
\ee
The maximum of the last factor, 
$\exp(-\eps^2/(2 \kappa))/\kappa^{\frac{3}{2}}$
is achieved for $\kappa=\eps^2/3$, and is valued
$\left( \frac{3}{e\, \eps^2} \right)^{\frac{3}{2}}$.

We still have a choice of $\eps=\Theta(1)$. This tuning is not crucial
for the complexity, however, it is easy to just choose the optimal
value.  As the saddle-point calculation predicts
$\kappa(1+\frac{A_0}{A_{\neq}}) n'= \eee(k+m) = n'/\bar{v}$, we have
\be
\frac{\eps^2}{3}
=
\kappa
=
\frac{A_{\neq}}{\bar{v}}
\ef,
\ee
and thus, for this value,
\be
\label{eq.53787536545d4}
\begin{split}
r_n(\{\tilde{u}_j\})
\lesssim
K_n
C
\frac{1}{n^{\frac{3}{2}}}
\frac{A_{\neq}}{\sqrt{2 \pi A_0}}
\left( \frac{\bar{v}}{e\, A_{\neq}} \right)^{\frac{3}{2}}
=
K_n
C
\frac{1}{n^{\frac{3}{2}}}
\sqrt{\frac{\bar{v}^3}{2 \pi e^3 \, A_0 A_{\neq}}}
\ef.
\end{split}
\ee
It is relatively easy to produce a value 
$K_n \sim C' n^{\frac{3}{2}}$ which certifies that
$r_n(\{\tilde{u}_j\}) \leq 1$ for all $\bt$, and that can be evaluated
efficiently. More details are given in Appendix \ref{app.contoKn}.
Using this value of $K_n$, the values $k$ and $m$ which are the
worst case, and also, up to small corrections, the ones predicted by
the saddle point, and the bound on the product of $\cosh$'s given by 
$\tilde{u}_j=\frac{n'}{k}$, gives
$r_n(\{\tilde{u}_j\}) \lesssim C''$ for some constant $C''<1$.

At this point the function $r_n$ has been chosen, and we are left with
the evaluation of $\eee(r_n(\{\tilde{u}_j\}))$. At this aim we do not
need anymore rigorous bounds, but the large-$n$ asymptotics is
sufficient to determine the asymptotic complexity.

The random values of $k$ and $m$ entering this average are
asymptotically Gaussian, with a quadratic form scaling with $n'$, that
is
$k-k^*, m-m^* \sim \sqrt{n'}$. In order to determine that the
algorithm is optimal up to a multiplicative constant, we only need to
know this fact. If we want to evaluate this constant, we have to
evaluate the quadratic form precisely. The quadratic form $Q_{(n,m)}$
associated to the fluctuations of $n$ and $m$ at fixed $k$, so that we
have (here and below the notation $X^*$, or $\eee(X)$, stands for the
saddle-point expectation of the random variable $X$)
\be
p_k(\delta n, \delta m)
\propto
\exp\left(
-\frac{1}{2k} 
(\delta n, \delta m) Q_{(n,m)}^{-1}
\begin{pmatrix}\delta n \\ \delta m \end{pmatrix}
\right)
\ee
is given by
\be
\begin{split}
Q_{(n,m)}
&=
\begin{pmatrix}
\eee((u_1+v_0 u_2)^2)
&
\eee((u_1+v_0 u_2)(u_2-v_0 u_1))
\\
\eee((u_1+v_0 u_2)(u_2-v_0 u_1))
&
\eee((u_2-v_0 u_1)^2)
\end{pmatrix}
\\
& \qquad
-
\begin{pmatrix}
\eee(u_1+v_0 u_2)^2
&
\eee(u_1+v_0 u_2)\eee(u_2-v_0 u_1)
\\
\eee(u_1+v_0 u_2)\eee(u_2-v_0 u_1)
&
\eee(u_2-v_0 u_1)^2
\end{pmatrix}
\end{split}
\ee
where expectations of $(u_1,u_2)=(v,\ell-1)$ are w.r.t.\ the measure
$\mu_{\neq}$. This is a direct consequence of the CLT.
The quadratic form 
$Q_{(k,m)}$
associated to the fluctuations of $k$ and $m$ at fixed $n$ are
obtained through a change of variables at leading order, described by
a matrix $B$
\begin{align}
Q_{(k,m)}^{-1}
&=
\frac{n}{k^*}
B^{\rm -T}
Q_{(n,m)}^{-1}
B^{\rm -1}
\ef;
&
B
&=
\begin{pmatrix} 
-\frac{k^*}{n} & 0 \\ -\frac{m^*}{n} & 1 
\end{pmatrix}
\ef.
\end{align}
That is, asymptotically,
\be
p_n(\delta k,\delta m)
\simeq
\frac{1}{2 \pi n \sqrt{\det Q_{(k,m)}}}
\exp\left(
-\frac{1}{2n} 
(\delta k, \delta m) Q_{(k,m)}^{-1}
\begin{pmatrix}\delta k \\ \delta m \end{pmatrix}
\right)
\ef.
\ee
The simple binomial factor
\be
\label{eq.86554264}
K_n \frac{(k+m-1)!}{k!\, m!} A_{\neq}^k A_0^m
\ee
also gives a Gaussian, with a quadratic form scaling with $n$, and
normalised so to be a constant of order 1 on the saddle-point values.
Again, this is enough for establishing the optimality of the algorithm,
but if we want to know the overall multiplicative constant, we have to
evaluate the associated quadratic form, that is obtained from Stirling
approximation:
\begin{align}
\label{eq.86554264b}
K_n \frac{(k+m-1)!}{k!\, m!} A_{\neq}^k A_0^m
&\sim
C''
\exp\left(
-\frac{1}{2n} 
(\delta k, \delta m) R
\begin{pmatrix}\delta k \\ \delta m \end{pmatrix}
\right)
\ef;
&
R
&=
\frac{n}{k^*}
\begin{pmatrix}
\frac{A_0}{A_{\neq}} & 1
 \\ 
1 & \frac{A_{\neq}}{A_0}
\end{pmatrix}
\ef.
\end{align}
Note that the matrix $R$ has rank 1, so it is only positive
semi-definite, with zero eigenvalue w.r.t.\ the infinitesimal
transformation 
$k/n \to k/n + \eps A_{\neq}$, 
$m/n \to m/n + \eps A_{0}$.

The most complicated correction factor is due to the product
\be
\frac{\prod_j 
\cosh \left(
\sqrt{\frac{3 A_{\neq}}{\bar{v}\,n'}}
\tilde{u}_j
\right)^{-1}
}
{
\cosh \left(
\sqrt{\frac{3 A_{\neq}}{\bar{v}\,n'}}
\frac{n'}{k^*}
\right)^{-k^*}
}
\ee
As said above, this quantity has been bounded by its evaluation at
$\tilde{u}_j=\eee(\tilde{u})=
\frac{n'}{k^*}=\frac{\bar{v}}{A_{\neq}}$, and the
overall constant coming from dealing with simplifications of the bound 
is the constant $C''$ appearing above in (\ref{eq.86554264b}).

The typical lists of $\tilde{u}_j$'s, however, 
have entries which are different, and follow some distribution with
all moments (and all cumulant moments) of order $1$, and all joint
moments almost factorised, i.e., for $j_{\alpha}$ distinct indices and
$\sum_{\alpha} \nu_{\alpha}$ of order 1,
\begin{align}
\eee(
\tilde{u}_{j_1}^{\nu_1}
\cdots
\tilde{u}_{j_h}^{\nu_h}
)
&=
M_{\nu_1} \cdots M_{\nu_h}
+\cO(n^{-1})
\ef;
&
M_1 &=
\frac{n'}{k^*}
=
\frac{\bar{v}}{A_{\neq}}
\ef.
\end{align}
Let us start our analysis with some heuristics.
As we have $k = \Theta(n)$, we are in a regime in which we can
approximate the functions $\cosh(\alpha/\sqrt{n})$ by
$\exp(-\alpha^2/(2n))$, which gives (recalling that $\eee(\tilde{u})k^*=n'$)
\be
\begin{split}
&
\eee
\left(
\frac{\prod_j 
\cosh \left(
\sqrt{\frac{3 A_{\neq}}{\bar{v}\,n'}}
\tilde{u}_j
\right)^{-1}
}
{
\cosh \left(
\sqrt{\frac{3 A_{\neq}}{\bar{v}\,n'}}
\frac{n'}{k^*}
\right)^{-k^*}
}
\right)
\simeq
\eee
\left(
\frac{
\prod_j 
\exp \left(
-\frac{3 A_{\neq}}{2\bar{v}\,n'}
\tilde{u}_j^2
\right)
}
{
\exp \left(
\frac{3 A_{\neq}}{2 \bar{v}\,n'}
\frac{{n'}^2}{(k^*)^2}
\right)^{-k^*}
}
\right)
\simeq
\frac{
\eee
\left(
\exp \left(
-\frac{3 A_{\neq}}{2\bar{v}}
\frac{1}{n'}
\sum_j
\tilde{u}_j^2
\right)
\right)
}
{
\exp \left(-\frac{3 A_{\neq}}{2 \bar{v}}
\frac{{n'}}{k^*}
\right)
}
\\
&
\qquad
\simeq
\frac{
\exp \left(
-\frac{3 A_{\neq}}{2\bar{v}}
\frac{1}{n'}
\eee
\left(
\sum_j
\tilde{u}_j^2
\right)
\right)
}
{
\exp \left(-\frac{3 A_{\neq}}{2 \bar{v}}
\frac{1}{n'}
\eee(\tilde{u})
\eee\left(\sum_j \tilde{u}_j \right)
\right)
}
=
\exp \bigg(
-\frac{3 A_{\neq}}{2\bar{v}}
\frac{1}{n'}
\eee
\bigg(
\sum_j
\big[
\tilde{u}_j^2
-
\eee(\tilde{u})
\tilde{u}_j
\big]
\bigg)
\bigg)
\\
& \qquad
\sim
\exp \bigg(
-\frac{3 A_{\neq}^2}{2\bar{v}^2}
\big[
\eee(\tilde{u}^2)
-
\eee(\tilde{u})^2
\big]
\bigg)
\ef.
\end{split}
\label{contoCHnonrig}
\ee
Some of the passages above, and especially the last one, are hard to
justify, as in principle we should consider possible correlations
between joint moments (that can be relevant even if they are of order
$1/n$, because, for example, in the expansion of an expression of the
form $\eee(\exp(\sum_{j=1}^n x_j))$ there are $\sim n^2$ moments of
the form $\eee(x_i x_j)$, and only $\sim n$ moments $\eee(x_i^2)$),
and correlations between moments of the $\tilde{u}_j$'s and the random
variable $k$.

However it is clear that these correlations may only affect the overall
constant, and cannot change the fact that the result is of order 1.

In order to perform a more rigorous calculation, we shall revert to
the Markov Chain approach that we have already used for the evaluation
of $\ppp(\textrm{reaching $D_n$})$ in (\ref{eq.ContoPDn}). Now we
rather approximate the functions $\cosh(\alpha/\sqrt{n})^{-1}$ by
$\cos(\alpha/\sqrt{n})$, and assume that, for a probability
distribution $p_a$,
\be
\eee_{\sum_j a_j = n}
\bigg(
\prod_j 
\cos\Big(\sqrt{\frac{\eps}{N}} a_j \Big)
\bigg) 
\sim
\exp(-n \xi_{N})
\ee
for some $\xi_{N}$ of order $1/N$. Self-consistency gives
\be
\exp(-n \xi_{N})
=
\sum_a
p_a
\exp(-(n-a) \xi_{N})
\cos\Big(\sqrt{\frac{\eps}{N}} a_j \Big)
\ee
that is
\be
1=
\frac{1}{2}
\left(
\hat{p}\Big(
\exp\Big(
\xi_{N}+i \sqrt{\smfrac{\eps}{N}}\Big)\Big)
+
\hat{p}\Big(
\exp\Big(
\xi_{N}-i \sqrt{\smfrac{\eps}{N}}\Big)\Big)
\right)
\ef.
\ee
This is consistent with the ansatz $\xi_N = \xi/N$, and the
result, obtained from Taylor expansion at
first non-trivial order,
\be
\xi=\frac{\eps^2}{2} \frac{\hat{p}'(1)+\hat{p}''(1)}{\hat{p}'(1)}
=\frac{\eps^2}{2} \frac{\eee_{p}(a^2)}{\eee_{p}(a)}
\ef.
\ee
Our seeked quantity at numerator is thus
\be
\exp(-n \xi_{n})
=
\exp(-\xi)
=
\exp\Big(
-\frac{\eps^2}{2} \frac{\eee_{p}(a^2)}{\eee_{p}(a)}
\Big)
\ee
while the quantity at denominator is, using $k^* \eee_{p}(a)=n'$,
\be
\exp\Big(
-\frac{\eps^2}{2 n'} k^* \eee_{p}(a)^2
\Big)
=
\exp\Big(
-\frac{\eps^2}{2} \eee_{p}(a)
\Big)
\ef.
\ee
Taking the ratio, and using the definition of $\eps^2$, gives back the
quantity obtained in (\ref{contoCHnonrig}).

As a result, the product of $\cosh$'s, at leading order, does not
interfere with the Gaussian average over $\delta k$ and $\delta m$,
and the overall acceptance ratio is given by
\be
\eee(r_n)
\simeq
\frac{C''}{\sqrt{\det (\bI + Q_{(k,m)} R)}}
\exp \bigg(
-\frac{3 A_{\neq}^2}{2\bar{v}^2}
\big[
\eee(\tilde{u}^2)
-
\eee(\tilde{u})^2
\big]
\bigg)
\ef.
\ee
The three factors above are all real values in $[0,1]$, and of order $1$.
This average rate, multiplied by the quantity
$\ppp(\textrm{reaching $D_n$})$ in (\ref{eq.ContoPDn}), gives the
inverse of the multiplicative factor in the complexity, w.r.t.\ the
Shannon bound.



                      \appendix

\section{A simple fact in Perron--Frobenius Theory}
\label{app.frobfact}

\noindent
This section is devoted to the proof of the following fact:
\begin{lemma}[Monotonicity of the Frobenius eigenvalue]
\label{lem.frobAB}
For $x \in \bR^+$ and $K=(\vec{a}|B)$ an irreducible non-negative
square matrix, let $K(x)$ be the matrix $K(x)=(x \vec{a}|B)$, and
$\lam_1(x)$ its Frobenius eigenvalue. The function $\lam_1(x)$ is
strictly increasing.
\end{lemma}

\noindent
{\it Proof.}  Call $\vec{v}_i(x)$, $\vec{v}_i(x)^{\ast}$ and
$\lam_i(x)$ the right-eigenvectors, left-eigenvectors and eigenvalues
of $K(x)$, respectively, with index $1$ reserved to the Frobenius one
(which is unique).  At first order in $\eps$ we have
\begin{align}
\vec{v}_1(x+\eps)
&=
\vec{v}_1(x)+\eps \sum_{j} \gamma_j(x) \vec{v}_j(x)
\ef;
&
\lam_1(x+\eps)
&=
\lam_1(x)+\eps \delta
\ef.
\end{align}
Decompose the vector $\vec{a}$ in the eigenbasis:
\begin{align}
\vec{a}
&=
\sum_{j} \alpha_j(x) \vec{v}_j(x)
\ef;
&
\alpha_j(x) 
&= \frac{(\vec{v}_j(x)^{\ast},\vec{a})}
{(\vec{v}_j(x)^{\ast},\vec{v}_j(x))}
\ef.
\end{align}
Note that, as $\vec{a}$ is positive and 
$\vec{v}_1(x)^{\ast}$ is Frobenius, $\alpha_1(x)>0$ for all $x>0$.
The defining equation
\be
K(x+\eps) 
\vec{v}_1(x+\eps)
=
\lam_1(x+\eps)
\vec{v}_1(x+\eps)
\ee
gives, at first order in $\eps$,
\be
\sum_j 
\lam_j(x) \gamma_j(x) \vec{v}_j(x)
+
(\vec{v}_1(x))_1
\sum_j 
\alpha_j(x) \vec{v}_j(x)
-
\delta
\vec{v}_1(x)
-
\lam_1(x)
\sum_j 
\gamma_j(x) \vec{v}_j(x)
=0
\ef.
\ee
Taking a scalar product with $\vec{v}_1(x)^{\ast}$ (and using the
unicity of the Frobenius vector for orthogonality) gives
\be
\delta
=
(\vec{v}_1(x))_1
\alpha_1(x)
>
0
\ef.
\ee
\hfill $\square$

\section{The {\sc BalancedShuffle} algorithm}
\label{app.BBHL}

\noindent
Here we describe an algorithm, given by Bacher, Bodini, Hollender and
Lumbroso in \cite[sec.~1.1]{axel2}, for the uniform sampling of
strings $\bnu=(\nu_1,\ldots,\nu_n) \in \{0,1\}^n$ with 
$\sum_i \nu_i=k$.
As stated in the text, this algorithm is `optimal' for the
randomness resource, in the sense that the random-bit complexity at
leading order coincides with the Shannon bound.

Call $\beta=k/n$.  The idea is that you sample the $\nu_i$'s one by
one, as if $\bnu$ were to be sampled with the measure 
$\mu(\bnu) = \bern_\beta^n$, (this costs $\Theta(1)$ per variable),
and, as soon as you have an excess of $\nu_i$'s equal to 0 or to 1,
complete with a sequence of Fisher--Yates random swaps. These swaps
cost $\sim \ln n$ each, but are performed on average only $\sim
\sqrt{n}$ times, so the swap part of the algorithm has subleading
complexity, and the overall complexity is genuinely linear, with no
logarithmic factors.  For completeness, we give in
Algorithm~\ref{algo.BBHL} a summary of the main features. In this
algorithm 
$\rint_{n}$ returns a uniform random integer in
$\{1,\ldots,n\}$.

\begin{algorithm}[H]
\Begin{
$a=k, \ b=n-k, \ i=0$\;
\Repeat{$a<0$ or $b<0$}
{
$i\,${\tt ++}\;
$\nu_i \longleftarrow \bern_{\beta}$\;
\lIf{$\nu_i=1$}{$a\,${\tt --}}\lElse{$b\,${\tt --}}
    }
\lIf{$a<0$}{$\bar{\nu}=0$}\lElse{$\bar{\nu}=1$}
\For{$j\gets i$ \KwTo $n$}{
$\nu_j=\bar{\nu}$\;
$h \longleftarrow \rint_{j}$\;
swap $\nu_j$ and $\nu_h$\;
    }
\Return{$\bnu$}
}
\caption{BBHL-shuffling.
\label{algo.BBHL}}
\end{algorithm}

\noindent
Note that, if $\beta$ is almost a $2$-adic number, namely
$\beta=\frac{a}{2^d}+ \eps$ for some integers $a$ and $d$, and
$\eps=o(n^{-\frac{1}{2}} (\ln n)^{-1})$, it may be convenient to just
use the value $\beta=\frac{a}{2^d}$, which speeds up the main part of
the algorithm, at the price of slowing down the subleading part.

\section{Certification of the constant $K_n$}
\label{app.contoKn}

\noindent
In Section \ref{sec.anAR}, equation (\ref{eq.53787536545d2}), we
need to produce an explicit expression $K_n$, as large as possible,
and such that we can certify
\be
\label{eq.5764752643x}
\begin{split}
K_n
\max_{k,m}
C
\frac{
(k+m-1)!}{k!\, m!}
A_{\neq}^k A_0^m
\cosh
\Big(
\frac{\eps \sqrt{n}}{k}
\Big)^{-k}
\leq 1
\ef.
\end{split}
\ee
(In this section, for simplicity of notation, we replace $n'$ by $n$,
so it is intended that $K_n$ here corresponds to $K_{n+v_0}$ in the
body of the paper).
Here $C$ and $\eps=\sqrt{3A_{\neq}/\bar{v}}$ are explicit constants,
and $A_0$, $A_{\neq}$ are parameters in $[0,1]$, with
$A_0+A_{\neq}=1$. All of these quantities are $\Theta(1)$, in a
scaling where $k,m,n$ are large.

In that section we anticipated the conclusion that, as needed for our
algorithm to reach optimal complexity up to a multiplicative constant,
we can choose $K_n$ as large as $\sim n^{\frac{3}{2}}$, which is
optimal up to a multiplicative constant as it is the scaling that this
quantity would have if we restricted the analysis to the sole position
of the saddle-point prediction, that is, calling 
$k=\kappa n$, $m=\mu n$ (with $\kappa, \mu = \Theta(1)$)
\begin{align}
\kappa
&=
\frac{A_{\neq}}{\bar{v}}
\ef;
&
\mu
&=
\frac{A_{0}}{\bar{v}}
\ef;
\end{align}
where, at leading order in a large-$n$ expansion, we have
(cf.\ equation(\ref{eq.53787536545d4}))
\be
\label{eq.53787536545d4app}
\begin{split}
K_n
=
n^{\frac{3}{2}}
C^{-1}
\sqrt{\frac
{2 \pi e^3 \, A_0 A_{\neq}}
{\bar{v}^3}
}
\ef.
\end{split}
\ee
This section is devoted to the description of a method by which we can produce a
constant $C'=\Theta(1)$, and a certification that, by choosing
$K_n \sim C' n^{\frac{3}{2}}$, equation (\ref{eq.5764752643x}) is satisfied.

First of all, while the original domain of maximisation is constituted
by $k,m$ integers, with $k \geq 1$ and $m \geq 0$, we can first extend
it to $m \in \bR^+$ (as this produces an upper bound, which is the
appropriate bound direction). In this case convexity in $m$, at given
$k$, is immediate, and we have a rather explicit equation for the
optimal value $m^*(k)$:
\be
\label{eq.spm}
\psi(m+k)-\psi(m+1)+\ln A_0 = 0
\ef,
\ee
where $\psi(x)=\deenne{x}{{}} \ln \Gamma(x)$ is the digamma function,
which, for $x$ large, goes as $\psi(x)=\ln x+\frac{1}{2x}+\ldots$, this
justifying one of the two saddle-point equations, $m/(k+m)=A_0$.
So we can rewrite (\ref{eq.5764752643x}) as the stronger condition
\be
\label{eq.5764752643xx}
\begin{split}
K_n
\max_{k \geq 1}
C
\frac{
\Gamma(k+m^*(k))}{k!\, \Gamma(m^*(k)+1)}
A_{\neq}^k A_0^{m^*(k)}
\cosh
\Big(
\frac{\eps \sqrt{n}}{k}
\Big)^{-k}
\leq 1
\ef.
\end{split}
\ee
This is an equation of the form $K_n \max_k g(k,m(k)) \leq 1$.
Now, suppose that we can certify that,
for some integers $a \leq b$,
$\dienne{k}{{}}g(k,m(k))>0$ for all $k \leq a$ and
and
$\dienne{k}{{}}g(k,m(k))<0$ for all $k \geq b$.
This implies that the maximum $k^*=k^*(n)$ is attained for $k \in
\{a,a+1,\ldots,b\}$.  As the evaluation of any single value of $g$, or
derivative, with $\cO(\ln n)$ digits of precision, takes a time of
order $\cO((\ln n)^{\gamma})$ for some $\gamma$, under the assumptions
above finding the maximum has complexity of order 
$\cO((\ln n)^{\gamma} (b-a+1))$, by just trying out all the values in
the potential interval of integers.

In general, we have
\begin{align}
\dienne{k}{{}} g(k,m(k))
&=
\deenne{k}{{}} g(k,m(k))
+
\deenne{m}{{}} g(k,m(k))
\,
\deennearg{k}{{}}{m(k)}
\ef.
\end{align}
However, in our application the function $m(k)$ is defined by the
equation $\deenne{m}{{}} g(k,m(k))=0$, so the expression above
simplifies into
\begin{align}
\dienne{k}{{}} g(k,m(k))
&=
\deenne{k}{{}} g(k,m(k))
\ef.
\end{align}
The quantity of interest, $\dienne{k}{{}} g(k,m(k))$, is thus,
up to overall positive factors,
\be
\big(
\psi(m^*(k)+k)-\psi(k+1)+\ln A_{\neq}
\big)
+
\left(
-\ln 
\cosh\Big(
\frac{\eps \sqrt{n}}{k}
\Big)
+
\frac{\eps \sqrt{n}}{k}
\tanh\Big(
\frac{\eps \sqrt{n}}{k}
\Big)
\right)
\ee
This quantity is of order $1/n$, as it can be seen by setting
$k=\kappa n$ and $m^*(k)=\kappa n A_0/A_{\neq} + \delta$, which
defines $\delta$ in terms of $\kappa$ through
(\ref{eq.spm}). At leading order in $1/n$, it changes sign where it is
expected to do so, i.e.\ at the position of the saddle point.
So, if we can bound the involved functions with absolute errors of
order $1/n^2$ (i.e., relative errors of order $1/n$), we
can have a bound on the position of the point(s) where the derivative
changes sign, which is of order $1/n$ on the scale of $\kappa$, that
is, of order $1$ on the scale of $k$, or, in other words,
$b-a+1=\Theta(1)$.

One involved function is $\psi(m+k)-\psi(m+1)$. This expression is
bounded at arbitrary order, by using (the analytic continuation of)
$\psi(x+N)-\psi(x)=\sum_{k=0}^{N-1} \frac{1}{x+k}$, and then the
Euler--Maclaurin formula at the desired order, with the customary
simplest bound on the error term,
\begin{align}
\sum_{i=m}^n f(i)
&
\lessgtr
\int_{m}^n \dx{x} f(x)
+
\frac{f(n)+f(m)}{2}
+
\sum_{k=1}^p
\frac{B_{2k}}{(2k)!}
\big(
f^{(2k-1)}(n)-f^{(2k-1)}(m)
\big)
\pm
|R_p|
\ef;
\\
|R_p|
&
:=
\frac{2 \zeta(2p)}{(2 \pi)^{2p}}
\int_{m}^n \dx{x} |f^{(2p)}(x)|
\ef.
\end{align}
This gives in particular, already at $p=1$,
\begin{align}
\psi(x+y+1)-\psi(x)
&
\lessgtr
\ln \frac{x+y}{x}
+
\frac{2x+y}{2x(x+y)}
+
\left(-\frac{1}{12}
\pm
\frac{2 \zeta(2p)}{(2 \pi)^{2p}}
\right)
\frac{y(2x+y)}{x^2(x+y)^2}
\ef;
\end{align}
which is of the desired form.
A second special function that appears in our estimates is
\be
c(x) := -\ln \cosh \sqrt{x} + \sqrt{x} \tanh \sqrt{x}
\ef.
\ee
This is a rather regular function, that behaves as $c(x)\sim x/2$ for
$x$ small and $c(x)\sim \ln 2$ for $x$ large. We can find, for example
\be
\frac{3 x (12 + x)}{(6 + x) (12 + 5 x)}
\leq
c(x)
\leq
\frac{3 x (6 + x)}{4 (3 + x)^2}
\ee
and the difference between the two bounds is of order $\sim x^4/864$,
which, in our case $x=\cO(1/n)$, scales as $n^{-4}$ and is far
sufficient at our aims.


\end{document}